\documentclass[12pt]{article}
\usepackage{amsmath,amssymb,amsfonts,amsthm}
\usepackage[utf8]{inputenc}
\usepackage[T1]{fontenc} 
\usepackage{enumerate}
\usepackage{graphicx}
\usepackage{tabularx}
\usepackage{rotating}
\usepackage{overpic}
\usepackage{a4wide}
\usepackage{color, xcolor}
\usepackage{url, comment, hyperref}

\numberwithin{equation}{section}

\def\R{\mathbb{R}}

\def\N{\mathbb{N}}
\def\Z{\mathbb{Z}}

\def\S{\mathbb{S}}

\newcommand{\MyParagraph}[1]{\medskip\noindent\textbf{#1}\;}
\newcommand{\dd}{\mathrm{d}}

\newcommand{\abs}[1]{\vert#1\vert}

\begin{document}
\title{\textbf{Community Integration Algorithms (CIAs)\\ for Dynamical Systems on Networks}}
\author{Tobias B{\"o}hle\footnote{Technical University of Munich, School of Computation Information and Technology, Department of Mathematics, Boltzmannstraße 3, 85748 Garching, Germany.}~,~
Mechthild Thalhammer\footnote{Leopold--Franzens Universit{\"a}t Innsbruck, Institut f{\"u}r Mathematik, Technikerstraße 13/7, 6020 Innsbruck, Austria.}~,~
Christian Kuehn\footnote{Technical University of Munich, School of Computation Information and Technology, Department of Mathematics, Boltzmannstraße 3, 85748 Garching, Germany.}
\footnote{Complexity Science Hub Vienna, Josefst{\"a}dterstraße 39, 1080 Vienna, Austria.}}
\date{\today}
\maketitle
\MyParagraph{Abstract.}
Dynamics of large-scale network processes underlies crucial phenomena ranging across all sciences. Forward simulation of large network models is often computationally prohibitive. Yet, most networks have intrinsic community structure. We exploit these communities and propose a fast simulation algorithm for network dynamics. In particular, aggregating the inputs a node receives constitutes the limiting factor in numerically simulating large-scale network dynamics. We develop community integration algorithms (CIAs) significantly reducing function-evaluations. We obtain a substantial reduction from polynomial to linear computational complexity. We illustrate our results in multiple applications including classical and higher-order Kuramoto-type systems for synchronisation and Cucker--Smale systems exhibiting flocking behaviour on synthetic as well as real-world networks. Numerical comparison and theoretical analysis confirm the robustness and efficiency of CIAs.
\section{Introduction}

The importance of, as well as the level of research activity in, network dynamics has seen a dramatic increase within the 21st century~\cite{WattsStrogatz,BarabasiAlbert}. Previously, network models had mostly focused on either (a) all-to-all coupling, or (b) on highly structured sparse cases such as classical lattices or trees. Simulating dynamics, e.g., synchronization~\cite{PikovskyRosenblumKurths}, collective motion~\cite{VicsekZafiris}, or contact processes~\cite{KissMillerSimon}, on networks of the type (a)-(b) is already non-trivial. Indeed, usually each node receives/collects inputs from its neighbours at each time step, processes this information, usually by some form of averaging, and then adjusts its own behaviour. From a computational viewpoint these steps are straightforward for very sparse interactions because there are just few function evaluations and the complexity scales linearly with the number of nodes. Already for all-to-all coupled networks, the processing step of averaging incoming information from neighbours is very costly. A direct computational approach yields that the number of function evaluations within the averaging step of processing the information at each node grows at least quadratically. This situation gets much worse for temporal networks~\cite{HolmeSaramaki}, adaptive/co-evolutionary networks~\cite{GrossSayama}, multiplex/multilayer~\cite{Boccalettietal} networks, and higher-order/polyadic interactions beyond graphs~\cite{Battistonetal1}. Yet, one might hope that there is a low-dimensional number of averaged order parameters (or observables), which are similar or even identical inputs for each node at each time step. 

Most real-world networks are neither extremely sparse nor extremely dense but rather contain many heterogeneous structures~\cite{Newman}. Therefore, using brute-force network simulations quickly encounters computational barriers. In this work, we are going to combine several mathematical ideas to simulate many large-scale network dynamics models efficiently. In this introduction, we start with a non-technical presentation of our approach involving two pre-simulation (``off-line'') steps (P1)-(P2) and two evaluation (``on-line'') steps (E1)-(E2). The more detailed technical development of the computational methodology starts in Section~\ref{sec:NetworkSystems}. The first step (P1) in our approach employs community detection algorithms to identify densely connected sub-networks. The second step (P2) is to approximate, if necessary, the coupling function between nodes via a common basis, e.g., using Fourier methods. This step helps us to identify possible observables. The order of (P1)-(P2) can be reversed or parallelized. Basically, (P1) tackles heterogeneity, while (P2) identifies the best observable to exploit local density within a community. For each community, we utilize the similarity of nodes to significantly reduce the information processing in step (E1) at each node, i.e., the \emph{local observable plays the role of a common input} reducing quadratic or worse polynomial scaling function evaluations to just linear cost within the number of nodes. Since our networks are assumed to be heterogeneous we also must account in step (E2) for the very sparsely connected nodes, which is possible by direct computation. Our approach yields significant reductions of the required memory capacities and the overall costs measured by the total numbers of function evaluations.

For a very specialized and particular case, we have demonstrated recently that employing simple variants of the steps (E1)-(E2) and (P1)-(P2) can work potentially work~\cite{BohleKuehnThalhammer}. In this work, we develop the general CIA framework and show that it works in an extremely broad class of network dynamics applications, that it is robustness with regard to real data sets, that the general method does yield linear computational complexity with respect to the dimensions of the systems, and that the steps naturally extend to higher-order/polyadic dynamics. The remaining parts are organised as follows. In Section~\ref{sec:NetworkSystems}, we introduce the considered classes of network dynamical systems. In Section~\ref{sec:CIA}, we detail and exemplify the key steps of Community Integration Algorithms (CIAs). In Section~\ref{sec:numerics}, we present the advantages of our approach and confirm the substantial gain in efficiency by a series of numerical experiments for widely-used models and real-world networks. This includes Kuramoto systems arising in the description of synchronisation, extended Kuramoto-type models involving higher-order/polyadic interactions, Cucker--Smale systems modelling collective motion, and collective motion on real-world animal networks. Generalisations to more complex frameworks and open questions are mentioned in Section~\ref{sec:conclusion}. Supplementary calculations and illustrations are collected in an appendix, which contains a detailed mathematical setup for all the examples as well as theoretical justification for the efficiency of CIAs.

\section{Dynamical Systems on Networks}
\label{sec:NetworkSystems}

Dynamical systems on networks are of importance in many sciences ranging from physics, chemistry, biology and medicine to social sciences~\cite{BarratBarthelemyVespignani,PorterGleeson}. Illustrative examples for time-continuous dynamical systems on networks include Desai-Zwanzig systems~\cite{DesaiZwanzig} describing the motion of interactive particles under the influence of external confining potentials, Kuramoto models~\cite{Kuramoto}, tracking the evolution of phase oscillators, Cucker-Smale systems~\cite{CuckerSmale1} describing the movements and flocking behavior of birds, coupled van-der-Pol/FitzHugh-Nagumo models frequently used in neuroscience~\cite{BoergersKopell}, and Hegelsenmann-Krause models for opinion formation~\cite{HegselmannKrause}. Instead of studying classical versions of these models using all-to-all coupling, we study several of these models on general networks that possess community structure. Even though the network models originate from different disciplines, they can all be described by one single general network model, which we focus on here. The model class is given by 
\begin{align}\label{eq:GeneralNetworkSystem}
    x_m'(t) = f_m(x_m(t)) + \frac 1N \sum_{\ell = 1}^N a_{m\ell}\ g(x_\ell(t),x_m(t)), \qquad '=\frac{\textnormal{d}}{\textnormal{d} t},~x_m(0)~\text{given},
\end{align}
where $m\in \{1,\dots,N\}$ and $t\in [0,T]$. This system is based on an underlying network that has $N$ nodes and is represented by an $N\times N$-dimensional adjacency matrix $A$ with entries $a_{m\ell}$. For the sake of simplicity, we restrict ourselves to an undirected and unweighted graph such that $A$ is additionally symmetric and $a_{m\ell}\in \{0,1\}$. Further, $x_m(t)$ denotes the state of node $m$ at time $t$, the functions $f_m$ describe the intrinsic dynamics of the $m$-th node and $g(x_\ell(t),x_m(t))$ is a general coupling function that describes the strength of the interaction that node $\ell$ has on node $m$, if they are coupled. Finally, $T>0$ denotes a final time until which we want to integrate the system \eqref{eq:GeneralNetworkSystem}. While the range of $x_m(t)$ and thus also the domain of $f_m$ and $g$ is generally part of an abstract space $\mathcal X$, we typically have $\mathcal X \in \{\R, \R^n, \R/(2\pi \Z), \dots \}$.
By grouping the states of the nodes into one common vector $x=(x_1,\dots,x_N)$ and introducing an $N$-dimensional function $H$ with components
\begin{subequations}
\label{eq:GeneralNetworkSystemH}
\begin{align}
    H_m(x) = f_m(x_m) + \frac{1}{N}\sum_{\ell=1}^N a_{m\ell}\ g(x_\ell, x_m), \quad m\in \{1,\dots,N\},
\end{align}
the initial value problem from \eqref{eq:GeneralNetworkSystem} can also be written as
\begin{align}
    x'(t) = H(x(t)).
\end{align}
\end{subequations}
Even though this is a very general formulation, many typical network models have special structure. For example, in many models, the coupling function $g$ is of the form $g(\tilde x,\hat x) = h(\tilde x - \hat x)$. Table \ref{tab:NetworkExamples}, that can be found in the appendix, shows all examples mentioned at the beginning of this section fit this framework. We remark that the classical version of an all-to-all coupling is retained as the special case of a complete network. This special case will be an automatically included sub-problem in our implementation of CIAs since the coupling within each community closely resembles an all-to-all coupling. Although the network systems that we have mentioned above are described by time-continuous dynamical systems and specifically by nonlinear ordinary differential equations, completely analogous considerations for CIAs hold for time-discrete network dynamics; we also cover the Bornholdt-Rohlf discrete-time network model for self-organized criticality~\cite{BornholdtRohlf} to illustrate this point. Furthermore, when numerically integrating a continuous-time dynamical system that is given by a system of ODEs, one first time-discretizes this ODE system. Our CIAs efficiently evaluate large sums that appear in the resulting time-discrete system. 

\section{Community Integration Algorithms}\label{sec:CIA}

To numerically integrate the system \eqref{eq:GeneralNetworkSystemH} one first discretizes the time interval $[0,T]$ into many small steps $0=t_0<t_1<\cdots<t_T = T$ and then employs an iterative time stepping scheme~\cite{BlanesCasas}, e.g., a Runge-Kutta or multistep method. Independent of the method, each time iteration step needs at least one evaluation of the right-hand side $H$. Therefore, it is of key importance for a fast numerical integration to implement the evaluation of $H$ efficiently. However, when one looks at the specific structure of $H$, one notices that each of its components consist of a large sum. In total there are $\mathcal O(N^2)$ operations (summations and evaluations of the coupling function $g$) necessary to evaluate $H(x)$ for a given $x$ only a single time. This quadratic dependence on $N$ severely restricts the number of nodes that a network can possess such that numerical simulations on it are tractable. For higher-order/polyadic systems, see e.g.~Appendix \ref{sec:HOKuramoto}, the situation even worsens.

\begin{figure}
    \centering
    \begin{overpic}[width=0.95\textwidth]{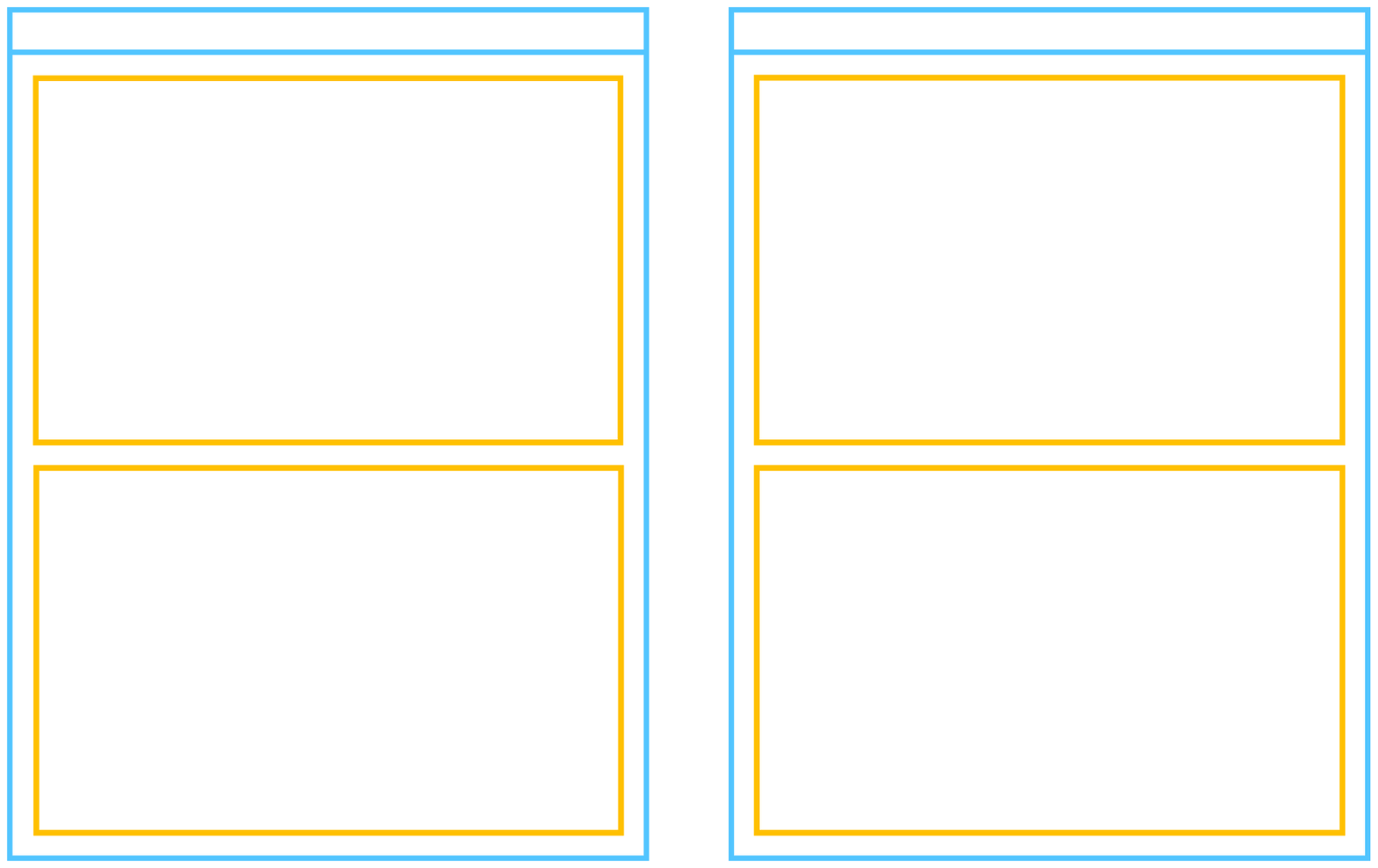}
    \put(11,60){\textbf{CIA Pre-simulation}}
    \put(68,60){\textbf{CIA Evaluation}}
    
    \put(10,54){Community Detection}
    \put(3,31){\includegraphics[width = 0.2\textwidth]{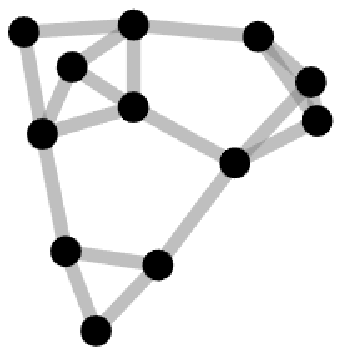}}
    \put(23.5,31){\includegraphics[width = 0.2\textwidth]{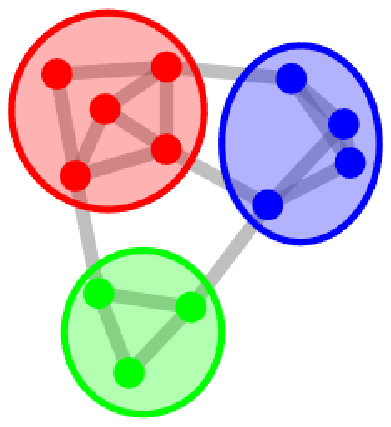}}
    \put(23,42){$\Rightarrow$}
    
    \put(7,26){High-Order Approximation}
    \put(3,2.5){\includegraphics[width=0.2\textwidth]{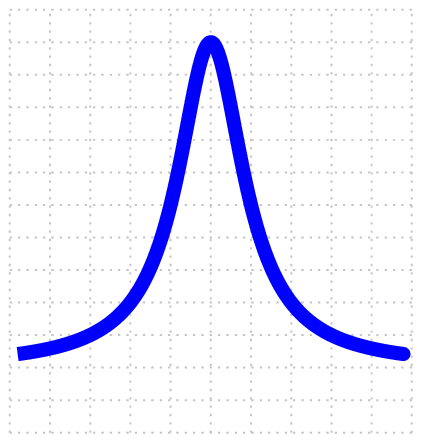}}
    \put(23.5,2.5){\includegraphics[width=0.2\textwidth]{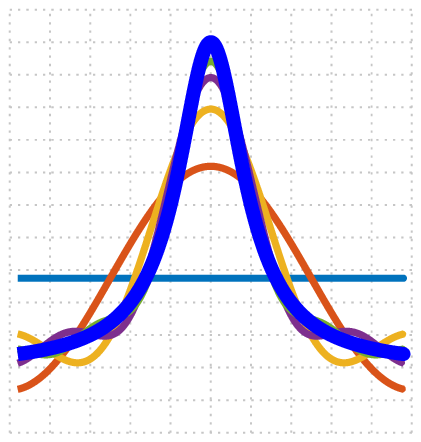}}
    \put(23,14){$\Rightarrow$}
    
    \put(56,54){Community Structure Exploitation}
    \put(66,31){\includegraphics[width = 0.2\textwidth]{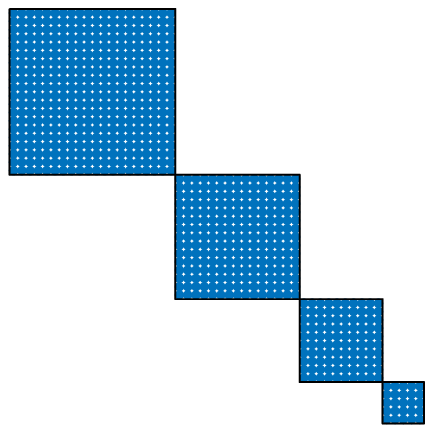}}
    
    \put(65,26){Sparse Summation}
    \put(66,2.5){\includegraphics[width = 0.2\textwidth]{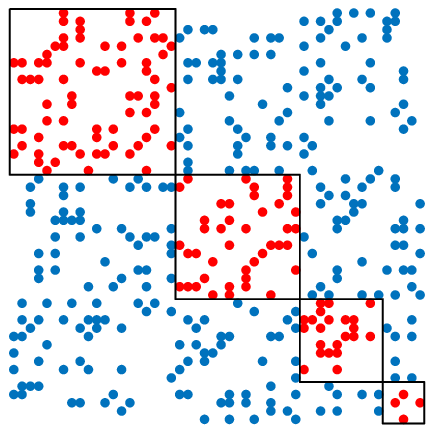}}

    \end{overpic}
    \caption{Flowchart that illustrates the steps of our Community Integration Algorithm.}
    \label{fig:flowchart}
\end{figure}

Our new Community Integration Algorithms (CIA) achieves to evaluate the right-hand side $H$ in \eqref{eq:GeneralNetworkSystemH} and requires only $\mathcal O(N)$ operations in each time step. It consists of four main steps. Two of them are done before the simulation and only need to be done once, whereas the other two have to be processed for each time step, see Figure \ref{fig:flowchart}. These steps are:

\begin{itemize}
\item[(P1)] Application of an effective community detection algorithm and transformation of the adjacency matrix by permutation to block form.
\item[(P2)] Identification of a suitable representation or high-order global approximation, respectively, of the coupling function $g$ to compute a suitable observable.
\item[(E1)] Exploiting community structure by computing a local observable for each community to avoid summations common among similar nodes.
\item[(E2)] Treatment of the remaining sparse parts of the network as well as small remaining heterogeneity within communities based on direct summations.
\end{itemize}

\begin{figure}
    \centering
    \begin{overpic}[width = \textwidth]{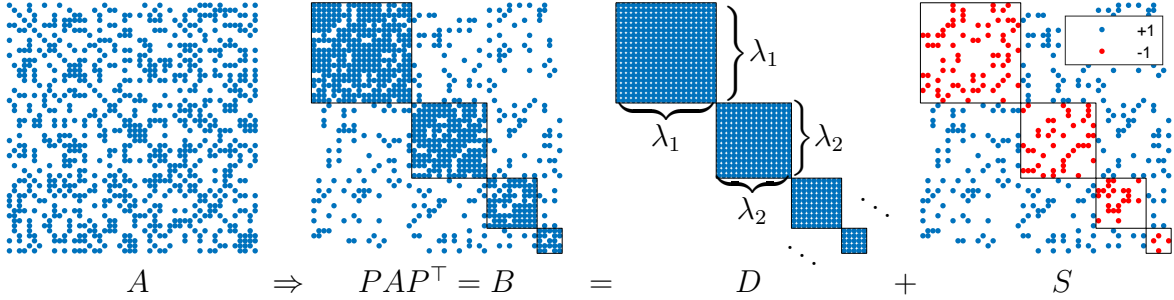}
    \put(12,3){$A$}
    \put(31,3){$PAP^\top=B$}
    \put(62,3){$D$}
    \put(88,3){$S$}
    
    \put(24,3){$\Rightarrow$}
    \put(50,3){$=$}
    \put(75,3){$+$}
    
    \put(61,22){$\Bigg\} \lambda_1$}
    \put(52.5,15){$\stackrel{\underbrace{\qquad\quad}}{\lambda_1}$}
    
    \put(66.5,15){$\bigg\} \lambda_2$}
    \put(60.5,9){$\stackrel{\underbrace{\qquad \; }}{\lambda_2}$}
    
    \put(72,9){$\displaystyle \ddots$}
    \put(66,5){$\displaystyle \ddots$}
    
    \end{overpic}
    \caption{Community detection step (P1)}
    \label{fig:Community_Detection}
\end{figure}

In the following we illustrate these steps via the example of phase oscillator systems
\begin{align}\label{eq:NetworkExample}
    \theta_m'(t) = f_m(\theta_m(t)) + \frac{1}{N} \sum_{\ell = 1}^N a_{m\ell}\  h(\theta_\ell(t)-\theta_m(t)).
\end{align}
Here $\theta_m(t)\in \S:=\R/(2\pi \Z)$ and the network is given by the adjacency matrix that can seen in Figure \ref{fig:Community_Detection}.

\subsection{CIA Pre-simulation Steps}
\MyParagraph{Community Detection (P1):} A key feature of our CIA is that it exploits the community structure of the underlying network. A community structure is present if the associated adjacency matrix is in block structure. However, when looking at the adjacency matrix that is depicted in Figure \ref{fig:Community_Detection}, there is at first no block structure evident. The community structure only becomes evident after permuting the nodes such that each community consists of nodes whose labels are consecutive integers. This permutation results in a matrix $B=PAP^\top$, where $P$ is a permutation matrix that is induced by a permutation $\kappa:\{1,\dots,N\}\to \{1,\dots,N\}$.

In many real-world scenarios the matrix $B$ does not have exact block structure but there may still be missing links inside a community and additional links across two communities. In any case, the matrix $B$ with evident block structure can be split into a dense matrix $D$ that has the exact block structure and a sparse matrix $S$. Here, $D$ only consists of entries in $\{0,1\}$ whereas $S$ is a sparse matrix with entries in $\{-1,0,1\}$. An $+1$ entry in $S$ denotes that there is an additional edge connecting two communities whereas a $-1$ entry represents a missing edge in a community.

The aim of step (P1) is to detect communities such that the permuted adjacency matrix $B$ has approximate block structure and can be decomposed into a dense matrix $D$ and a sparse matrix $S$, as seen in Figure \ref{fig:Community_Detection}. It is not necessary to store the full matrix $D$ but only the sizes of the communities that we denote by $\lambda_1, \lambda_2,\dots,\lambda_M$ if there are a total of $M$ communities. For numerical reasons it is sometimes better to allow nodes to belong to no community instead of letting them be part of a community that consists of very few or even only one node. Further, $S$ can be stored in a sparse format, so the total memory requirement is $\mathcal O(N)$. There are many effective algorithms that can detect communities in the underlying graph and thus transform the associated adjacency matrix to block structure, see Appendix \ref{sec:CommunityDetection}.

By using $PAP^\top=B=D+S$, where $P$ is the permutation matrix induced by the permutation $\kappa$, and $\phi_m(t) := \theta_{\kappa(m)}(t)$, we can write \eqref{eq:NetworkExample} as
\begin{align*}
    \theta_{\kappa(m)}'(t) &= f_{\kappa(m)}(\theta_{\kappa(m)}(t)) + \frac{1}{N}\sum_{\ell = 1}^N a_{\kappa(m)\kappa(\ell)}\ h(\theta_{\kappa(\ell)}(t)-\theta_{\kappa(m)}(t))
\end{align*}
    and so
\begin{align*}
    \phi_m'(t) &= f_{\kappa(m)}(\phi_m(t)) + \frac{1}{N}\sum_{\ell = 1}^N b_{m\ell}\ h(\phi_\ell(t)-\phi_m(t))\\
    &= \underbrace{f_{\kappa(m)}(\phi_m(t)) + \frac{1}{N}\sum_{\ell = 1}^N s_{m\ell}\ h(\phi_\ell(t)-\phi_m(t))}_{=:H^\text{sparse}_m(\phi(t))}   +    \underbrace{\frac{1}{N}\sum_{\ell = 1}^N d_{m\ell}\ h(\phi_\ell(t)-\phi_m(t))}_{=:H^\text{dense}_m(\phi(t))}.
\end{align*}
Thus, when evaluating the right-hand side, we only need to compute
\begin{align}\label{eq:NetworkExamplePhi}
    \phi_m'(t) = H_m^\text{sparse}(\phi(t)) + H_m^\text{dense}(\phi(t)).
\end{align}
Since $d_{m\ell} = 1$ if $1\le m,\ell\le \lambda_1$ and $d_{m\ell} = 0$ if $1\le m \le \lambda_1$ and $\ell > \lambda_1$ we obtain
\begin{align*}
    H^\text{dense}_m(\phi) = \frac{1}{N}\sum_{\ell = 1}^{\lambda_1} h(\phi_\ell-\phi_m)
\end{align*}
for $1\le m \le \lambda_1$, making it effectively an all-to-all coupling within that community. Similar representations of $H^\text{dense}_m(\phi)$ apply when $m>\lambda_1$.

\MyParagraph{High-Order Approximation (P2):}
This pre-simulation step is all about the expansion of the coupling function $h$. Since the domain is $2\pi$-periodic it makes sense to identify Fourier coefficients $a_k,b_k$, for $k\in \N$ with
\begin{align*}
    h(\phi) = \sum_{k=0}^\infty \Big( a_k \cos(k \phi) + b_k \sin(k\phi) \Big).
\end{align*}
Equivalently, one can also consider a Fourier expansion based on a complex Fourier series, see Appendix \ref{sec:Kuramoto}. Alternatively, yet not suitable here, an expansion in terms of polynomials is possible.
In any case for numerical reasons we terminate the series at a finite $k$ and only deal with the approximation
\begin{align*}
    h(\phi) \approx \sum_{k=0}^p \Big( a_k \cos(k \phi) + b_k \sin(k\phi)\Big),
\end{align*}
for some $p\in \N$, which determines the accuracy of the approximation. We shall see below that approximating the coupling function makes it easier to identify good local observables.

\subsection{CIA Evaluation Steps}

Based on the preparation done in the CIA pre-simulation steps, the right-hand side of the initial value problem \eqref{eq:NetworkExample} or equivalently \eqref{eq:NetworkExamplePhi} can now be evaluated using just $\mathcal O(N)$ operations. This evaluation is structured into two main steps:

\MyParagraph{Community Structure Exploitation (E1):}
This step aims to compute $H^\text{dense}(\phi)$, which, as shown previously, can be written as
\begin{align*}
    H^\text{dense}_m(\phi) = \frac 1N \sum_{\ell = 1}^{\lambda_1} h(\phi_\ell-\phi_m),
\end{align*}
if $m$ is a node belonging to the first community, i.e. $1\le m \le \lambda_1$ and similar representation are possible when $m$ is not in the first community. Combining this with the Fourier expansion that we conducted in (P2) and by using addition theorems for $\sin$ and $\cos$, we obtain
\begin{align}
    \nonumber
    H^\text{dense}_m(\phi) &\approx \frac{1}{N} \sum_{\ell = 1}^{\lambda_1} \sum_{k=0}^p \Big( a_k \sin(k(\phi_\ell - \phi_m)) + b_k \cos(k(\phi_\ell - \phi_m))\Big)\\
    \nonumber
    &=\sum_{k=0}^p \frac{1}{N} \sum_{\ell = 1}^{\lambda_1} \Big( a_k \sin(k\phi_\ell)\cos(k\phi_m)-a_k\cos(k\phi_\ell)\sin(k\phi_m)\\
    \label{eq:HdenseSeparate}
    &\qquad+ b_k \cos(k\phi_\ell) \cos(k\phi_m)+ b_k\sin(k\phi_\ell)\sin(k \phi_m)\Big),
\end{align}
for $1\le m \le \lambda_1$. Even though it first seems a lot more messy, we have separated terms involving $\phi_\ell$ and terms with $\phi_m$. Since we sum over $\ell$ and the terms involving $\phi_\ell$ remain the same for each $m$, we can precompute quantities
\begin{align}\label{eq:PrecomputeQuantities}
    q^\text{cos}_k := \frac 1N \sum_{\ell = 1}^{\lambda_1} \cos(k\phi_\ell) \quad \text{and} \quad q^\text{sin}_k := \frac 1N \sum_{\ell = 1}^{\lambda_1} \sin(k\phi_\ell),
\end{align}
for $k=0,\dots,p$. In particular, \eqref{eq:PrecomputeQuantities} are precisely the local observables, which are felt by all nodes within one community. They allow us to conclude
\begin{align}
    \nonumber
    H^\text{dense}_m(\phi) &\approx \sum_{k=0}^p\Big( a_k q^\text{sin}_k \cos(k\phi_m)-a_k q^\text{cos}_k\sin(k\phi_m)\\
    \label{eq:FastHdense}
    &\qquad+b_k q^\text{cos}_k \cos(k\phi_m) + b_k q^\text{sin}_k \sin(k\phi_m)\Big)
\end{align}
for $1\le m\le \lambda_1$. It is important to note that the computational complexity in this representation of $H_m^\text{dense}(\phi)$ in independent of the total number of oscillators $N$. In summary, the procedure in this step is as follows:
\begin{enumerate}
    \item\label{item:precompute} Precompute the quantities \eqref{eq:PrecomputeQuantities} for each $k=0,\dots,p$ and similar quantities for other communities.
    \item\label{item:evaluate} Use the precomputed quantities to obtain a high-order approximation of $H^\text{dense}(\phi)$ according to formula \eqref{eq:FastHdense} and equivalent formulas for other communities.
\end{enumerate}

The computational complexity of Step \ref{item:precompute} is $\mathcal O (\lambda_1 p)$ for the first community and consequently $\mathcal O(N p)$ for the whole step. The same applies to Step \ref{item:evaluate}.

\MyParagraph{Sparse Summation (E2):}
In this step first $H^\text{sparse}(\phi)$ is evaluated and then combined with the results from the previous step to obtain the final right-hand side of \eqref{eq:NetworkExamplePhi}.
Recall that $H^\text{sparse}_m(\phi)$ is given by
\begin{align*}
    H^\text{sparse}_m(\phi) = f_{\kappa(m)}(\phi_m(t)) + \frac{1}{N}\sum_{\ell = 1}^N s_{m\ell}\ h(\phi_\ell(t)-\phi_m(t)).
\end{align*}
Since the matrix $S$ with entries $s_{m\ell}$ is sparse with at most $\mathcal O(N)$ non-zero entries, it is evident that even a straightforward summation of $H^\text{sparse}(\phi)$ requires only $\mathcal O(N)$ operations. This is exactly what this step is supposed to do. Finally, in this step we compute $H^\text{dense}(\phi) + H^\text{sparse}(\phi)$, which gives the right-hand side of \eqref{eq:NetworkExamplePhi}. All of these computations can be done with a complexity of $\mathcal O(N)$.

\subsection{Extendibility to Other Network Models}

While the previous subsection only illustrate the CIA steps for one particular model \eqref{eq:NetworkExample}, it is straightforward to see that it is applicable to many more network models. Our main argument to support this claim is that by using a Fourier or polynomial expansion of the coupling function $g$ or $h$, the parts containing $\phi_\ell$ and those comprising $\phi_m$ or $x_\ell$ and $x_m$, when dealing with a model that does not have a circular domain, respectively, can always be separated, as done in \eqref{eq:HdenseSeparate}. This allows the precomputation of quantities that do not depend on $m$ but appear in each component of the right-hand side. These quantities consist of large sums whose single precomputation prevents unnecessary sums in the evaluation of each component of the right-hand side. Even when the coupling function $g$ is not of the form $g(\tilde x, \hat x) = h(\tilde x - \hat x)$, a two dimensional Fourier or polynomial expansion of $g$ in terms of $\tilde x$ and $\hat x$ is possible, see Appendix \ref{sec:GeneralNetworkFourier} and \ref{sec:GeneralNetworkPolynom}.
While a Fourier expansion is suitable for any kind of periodic domain, the choice between a Fourier and a polynomial expansion relies on the structure of the coupling function. If the coupling function is already an ordinary or trigonometric polynomial the choice is simple, see for example Appendix \ref{sec:DesaiZwanzig}. If that is not the case both methods generally make sense. An approximation of the coupling function $g$ then has to be chosen in a way that the approximation accurately resembles $g$ at the points where it is evaluated. If, for example, all the particles $x_m(t)$ stay in a subset of the domain, it is only necessary to approximate the coupling function there. Furthermore, the coupling function does not always need to be approximated by either a Fourier series or a polynomial. Rather, parts or components of $g$ that resemble a polynomial structure can be treated with a polynomial expansion while other parts and components might be better approximated by a Fourier series, see for example Appendix \ref{sec:CuckerSmale}. This further enlarges the class of network systems that can be integrated using out CIA method. 
Finally, we want to highlight that our method is also applicable to higher-order/polyadic interactions in which the coupling function $g$ depends on more than two arguments, see Appendix \ref{sec:HOKuramoto}.

\section{Numerical Results}\label{sec:numerics}

\begin{figure}[t]
    \centering
    \begin{overpic}[width=0.7\textwidth]{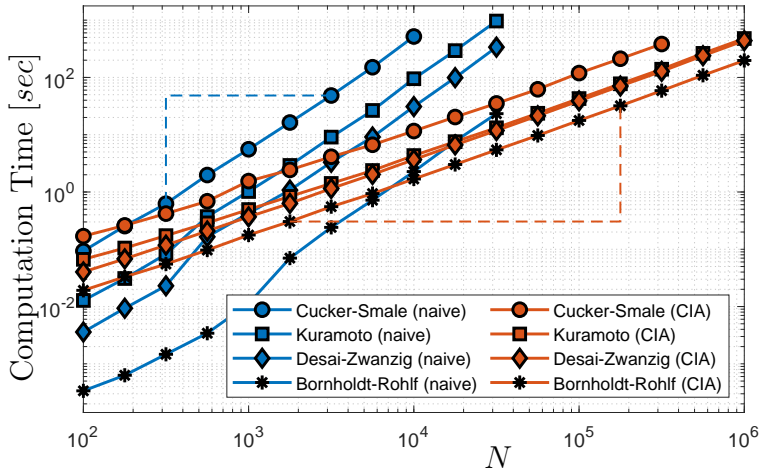}
        \put(60,0){$N$}
        \put(7,10){\begin{rotate}{90}Computation Time $[sec]$\end{rotate}}
    \end{overpic}
    \caption{Integration times for Cucker-Smale, Kuramoto, Desai-Zwanzig and Bornholdt-Rholf systems on a sequence of synthetically generated graphs of the form as shown in Figure \ref{fig:Community_Detection}. To integrate these systems, we used an explicit Euler scheme on an equidistant discretization of $[0,T]$ with $T=20$ and $\Delta t = 0.1$. The computation of the Cucker-Smale model takes the longest, since it has a two-dimensional phase space for each node. For $N\gtrsim 10^{4.5}$ storing the full matrix $A\in \R^{N\times N}$ exceeds memory capacities on current standard desktops. In this case, naive algorithms are not applicable. Caused by their lower memory requirement, CIAs can still cope with much larger graphs on standard desktop architectures.}
    \label{fig:model_comparison}
\end{figure}

Based on numerical simulations for a wide variety of widely-used large-scale network models we demonstrate below the advantages of using a CIA in comparison to a naive approach. Including synthetic and real-world networks and for pairwise and higher-order coupling we provide numerical evidence for robustness and efficiency of our approach. The network models that we consider in this section are
\begin{itemize}
    \item a Cucker-Smale model describing animal movement,
    \item a Kuramoto model for phases of oscillators on the unit circle,
    \item a Desai-Zwanzig model for interacting particles,
    \item and a Bornhold-Rholf model for self-organized criticality.
\end{itemize}
We compare the models on a fixed computational architecture using a sequence of networks that consist of four communities as seen in Figure \ref{fig:Community_Detection}. On the one hand, as seen in Figure \ref{fig:model_comparison}, the computation time for the naive approach depends quadratically on $N$ for all network models. Importantly, on the other hand, when using a CIA, the computation time depends only linearly on $N$. Furthermore, the memory requirements of a CIA are much lower as we can take advantage of sparsity outside of communities, while just calculating and storing a few observable values within each community, so much larger network sizes $N$ are possible;~cf.~Figure~\ref{fig:model_comparison}.

\begin{figure}[t]
    \centering
    \begin{overpic}[width = .4\textwidth]{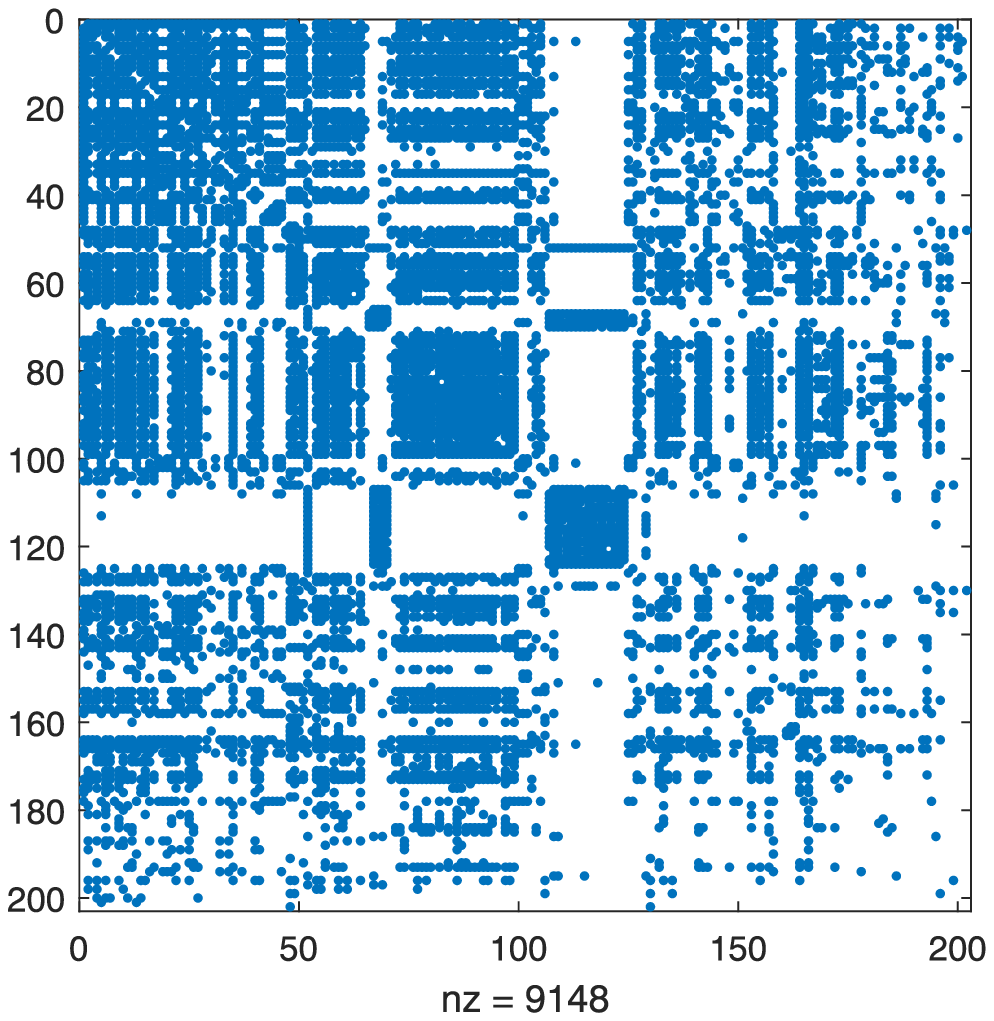}
        \put(7,0){\textbf{(a)}}
    \end{overpic}
    \begin{overpic}[width = .4\textwidth]{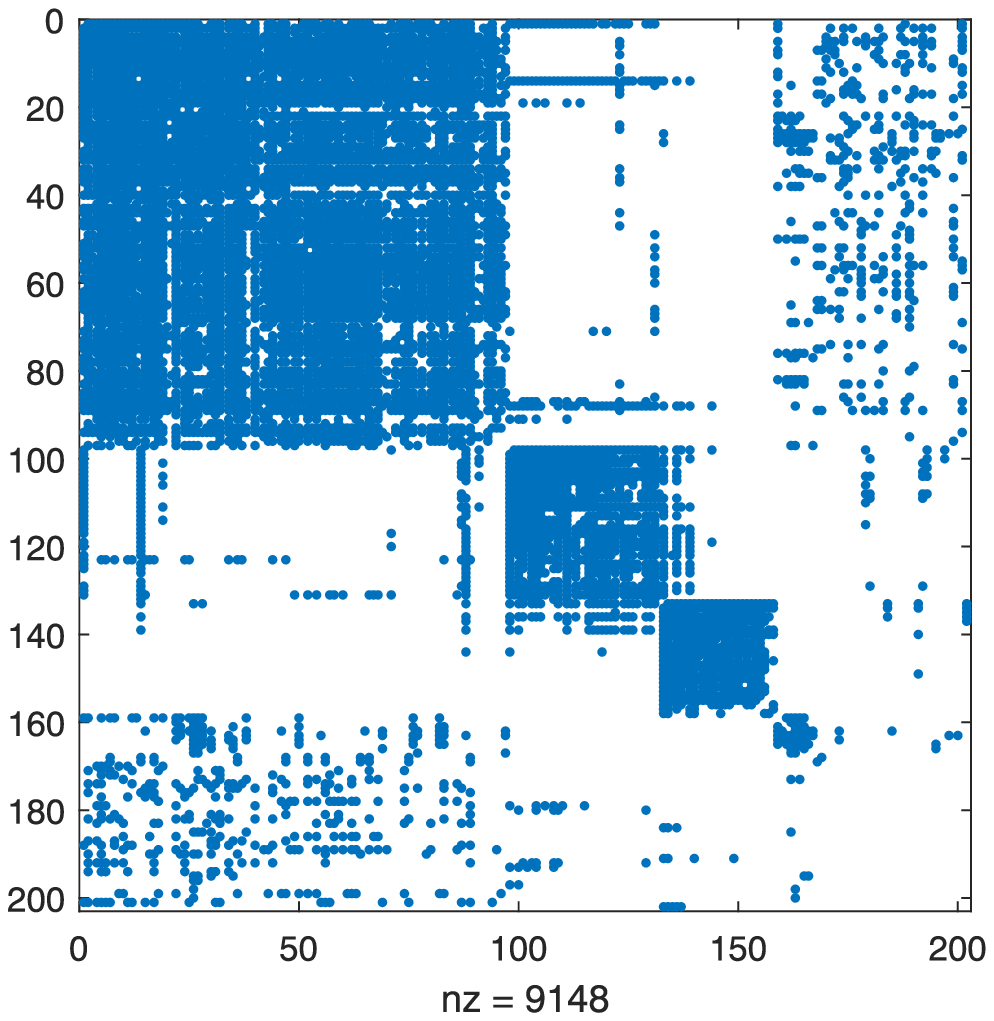}
        \put(7,0){\textbf{(b)}}
    \end{overpic}
    \caption{A real-world network from bird interactions \cite{Adelman2015, NetworkRepository}. \textbf{(a)} The adjacency matrix $A$ in its original form. \textbf{(b)} The adjacency after permuting the nodes such that a block structure is apparent.}
    \label{fig:RealWorldNetwork}
\end{figure}

\begin{figure}[t]
    \centering
    \begin{overpic}[width = .5\textwidth]{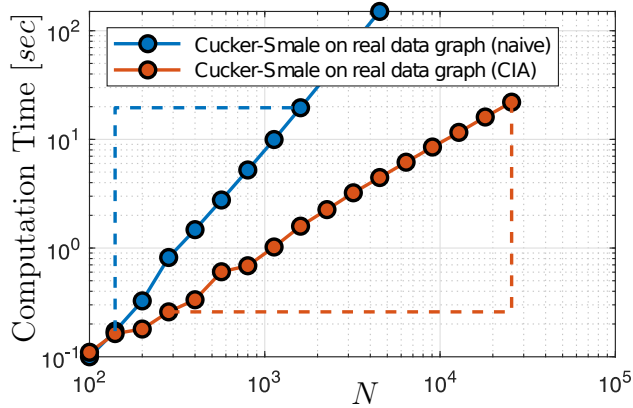}
        \put(52,0){$N$}
        \put(0,10){\begin{rotate}{90}Computation Time $[sec]$\end{rotate}}
    \end{overpic}
    \caption{A simulation of the Cucker-Smale model on a real-world network depicted in Figure \ref{fig:RealWorldNetwork}. To integrate these systems, we used an explicit Euler scheme on an equidistant discretization of $[0,T]$ with $T=20$ and $\Delta t = 0.1$.}
    \label{fig:CSTime_RealWorld}
\end{figure}

Next, we simulate the Cucker-Smale model on a real-world network that contains data from real bird interactions \cite{Adelman2015, NetworkRepository}, see Figure \ref{fig:RealWorldNetwork}. Based on this network, we construct a sequence of networks with growing sizes such that each network in the sequence still reflects the community structure of the original network and the amount of edges deviating from this community structure grows at most linearly in $N$. This helps us to study the effect of the network size on the computation time.
Again, as seen in Figure \ref{fig:CSTime_RealWorld}, the computation time of a naive approach of evaluating the right-hand side scales with $N^2$. The computational complexity of a CIA approach is only $\mathcal O(N)$. This clearly shows that a CIA can be be used to simulate dynamics on real-world networks.

Finally, we want to demonstrate that the idea of a CIA can also be applied to models that are beyond the general formulation \eqref{eq:GeneralNetworkSystem}. When considering higher-order models, such as a higher-order Kuramoto model, the computational complexity of a naive approach can be much worse than just $\mathcal O(N^2)$. However, even then a CIA approach is possible. One can still pre-compute suitable observables within communities and then evaluate the right-hand side based on these. This reduces the computational complexity to just $\mathcal O(N)$, see Figure \ref{fig:HOKuramoto}. For the details, see Appendix \ref{sec:HOKuramoto}. 

\begin{figure}[t]
    \centering
    \begin{overpic}[width = 0.5\textwidth]{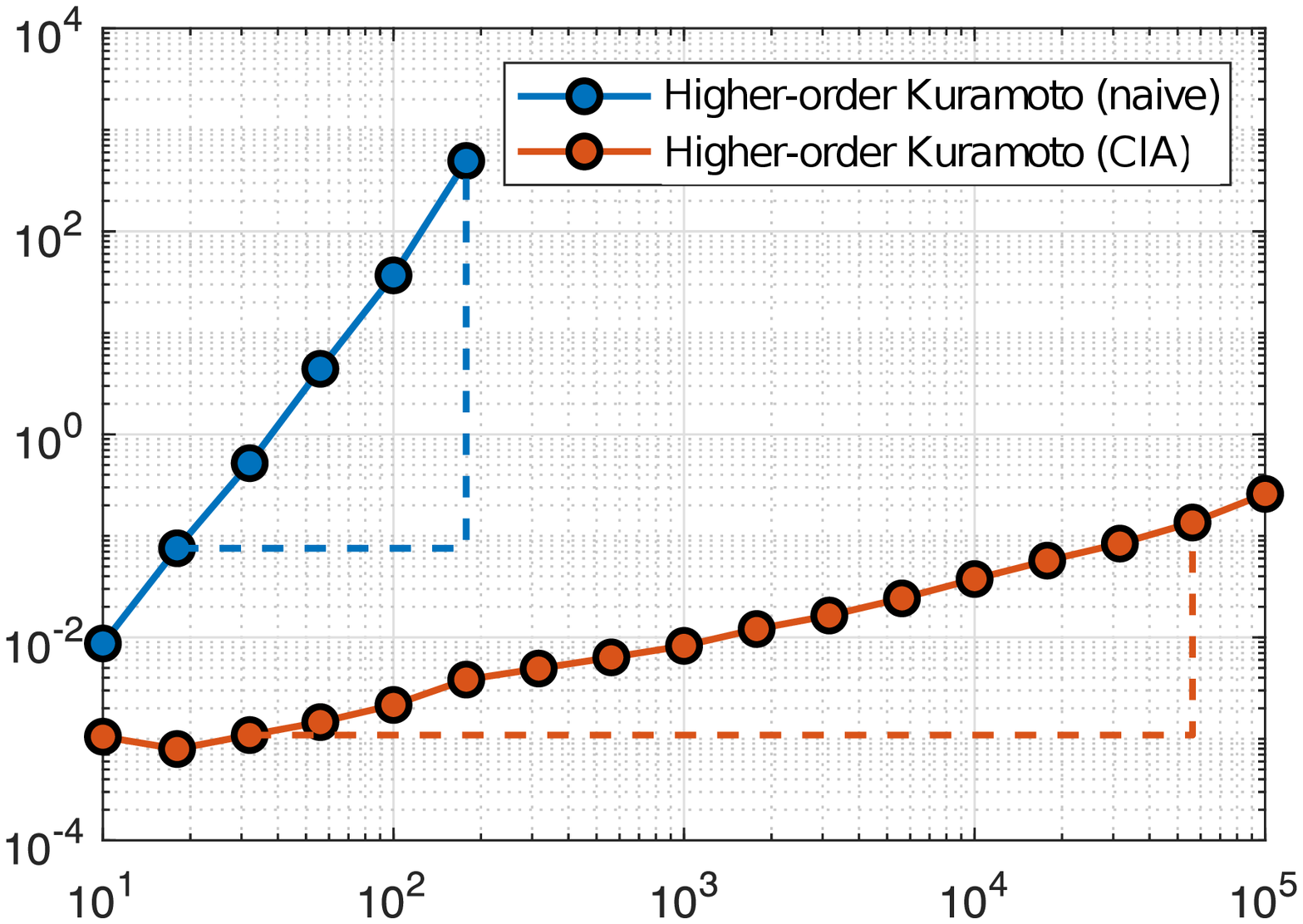}
        \put(58,0){$N$}
        \put(0,10){\begin{rotate}{90}Computation Time $[sec]$\end{rotate}}
    \end{overpic}
    \caption{Numerical integration of a higher-order Kuramoto system. A naive formulation and implementation leads to quartic complexity $\mathcal O(N^4)$ (blue), whereas the application of a community integration algorithm based on a suitable reformulation and precomputations permits the reduction to linear complexity $\mathcal O(N)$ (red).}
    \label{fig:HOKuramoto}
\end{figure}

\section{Conclusion and Outlook}
\label{sec:conclusion}

In summary, we have developed a new method to integrate large-scale network models with community structure. Unlike a naive approach, the computational complexity of CIAs is always linear in the number of involved nodes $N$. This can be achieved by first detecting a community structure in the underlying network and expanding the coupling function to identify local observables. Then, in each time step, based on precomputations of these observables and a reformulation of the right-hand side, it can be evaluated in just linear complexity $\mathcal O(N)$. As we have shown this approach is applicable to a wide variety of networks. Further, the approach works independent of whether the coupling is pairwise or higher-order, and whether the underlying network is synthetic or a real-world network, leading to myriad application across all areas of science.  

Lastly, we mention a few possible extensions: In the community detection step, we have so far focused of finding one community structure based on which we then perform the CIA evaluation steps (E1) and (E2). However, if the adjacency matrix has either additional and/or further hidden structure, additional computational savings are conceivable, see for example Appendix \ref{sec:ExtensionE1}. One can also generalize CIAs to networks whose adjacency matrices have large blocks away from the diagonal. For synthetically generated adjacency matrices in which this structure is already evident, our CIA steps work as well, but in real-world networks this structure first needs to be detected. Extensions are also possible for the high-order approximation step (P2). We considered approximations that are based on Fourier expansions, polynomials or a combination thereof. There might be even more functional approximation systems such that a separation of terms involving $x_m$ and $x_\ell$ in (E1) and thus a fast evaluation is possible.
\medskip

\noindent \textbf{Acknowledgements:} TB thanks the TUM Institute for Advanced Study (TUM-IAS) for support through a Hans Fischer Fellowship awarded to Christian Bick. TB also acknowledges support of the TUM TopMath elite study program. CK thanks the VolkswagenStiftung for support via a Lichtenberg Professorship.

\bibliographystyle{plain}
\bibliography{main_Tobi}

\begin{thebibliography}{10}

\bibitem{Adelman2015}
J.~S. Adelman, S.~C. Moyers, D.~R. Farine, and D.~M. Hawley.
\newblock {Feeder use predicts both acquisition and transmission of a
  contagious pathogen in a North American songbird}.
\newblock {\em Proceedings of the Royal Society B: Biological Sciences},
  282(1815), sep 2015.

\bibitem{BarabasiAlbert}
A.L. Barab{\'{a}}si and R.~Albert.
\newblock Emergence of scaling in random networks.
\newblock {\em Science}, 286(5439):509--512, 1999.

\bibitem{BarratBarthelemyVespignani}
A.~Barrat, M.~Barth\'{e}lemy, and A.~Vespignani.
\newblock {\em Dynamical Processes on Complex Networks}.
\newblock CUP, 2008.

\bibitem{Battistonetal1}
F.~Battiston, G.~Cencetti, I.~Iacopini, V.~Latora, M.~Lucas, A.~Patania, J.-G.
  Young, and G.~Petri.
\newblock Networks beyond pairwise interactions: structure and dynamics.
\newblock {\em Phys. Rep.}, 874:1--92, 2020.

\bibitem{BlanesCasas}
S.~Blanes and F.~Casas.
\newblock {\em A Concise Introduction to Geometric Numerical Integration}.
\newblock CRC Press, 2016.

\bibitem{Blondel2008}
V.~D. Blondel, J-L. Guillaume, R.~Lambiotte, and E.~Lefebvre.
\newblock {Fast unfolding of communities in large networks}.
\newblock {\em Journal of Statistical Mechanics: Theory and Experiment},
  2008(10):P10008, oct 2008.

\bibitem{Boccalettietal}
S.~Boccaletti, G.~Bianconi, R.~Criado, C.I. del Genio, J.~G{\'o}mez-Gardenes,
  M.~Romance, I.~Sendina-Nadal, Z.~Wang, and M.~Zanin.
\newblock The structure and dynamics of multilayer networks.
\newblock {\em Phys. Rep.}, 544(1):1--122, 2014.

\bibitem{BohleKuehnThalhammer}
T.~B{\"{o}}hle, C.~Kuehn, and M.~Thalhammer.
\newblock {On the reliable and efficient numerical integration of the Kuramoto
  model and related dynamical systems on graphs}.
\newblock {\em International Journal of Computer Mathematics}, pages 1--27, jul
  2021.

\bibitem{BoergersKopell}
C.~B{\"{o}}rgers and N.~Kopell.
\newblock Synchronization in networks of excitatory and inhibitory neurons with
  sparse, random connectivity.
\newblock {\em Neural Comput.}, 15(3):509--538, 2003.

\bibitem{BornholdtRohlf}
S.~Bornholdt and T.~Rohlf.
\newblock Topological evolution of dynamical networks: global criticality from
  local dynamics.
\newblock {\em Phys. Rev. Lett.}, 84(26):6114--6117, 2000.

\bibitem{Clauset2004}
A.~Clauset, M.~E.~J. Newman, and C.~Moore.
\newblock {Finding community structure in very large networks}.
\newblock {\em Physical Review E}, 70(6):066111, dec 2004.

\bibitem{CuckerSmale1}
F.~Cucker and S.~Smale.
\newblock Emergent behavior in flocks.
\newblock {\em IEEE Transactions on Automatic Control}, 52(5):852--862, 2007.

\bibitem{DesaiZwanzig}
R.C. Desai and R.~Zwanzig.
\newblock Statistical mechanics of a nonlinear stochastic model.
\newblock {\em J. Stat. Phys.}, 19:1--24, 1978.

\bibitem{GrossSayama}
T.~Gross and H.~Sayama, editors.
\newblock {\em Adaptive Networks: Theory, Models and Applications}.
\newblock Springer, 2009.

\bibitem{Hagberg2008}
A~A Hagberg, D~A Schult, and P~J Swart.
\newblock {Exploring network structure, dynamics, and function using NetworkX}.
\newblock {\em 7th Python in Science Conference (SciPy 2008)},
  (February):11--15, 2008.

\bibitem{HegselmannKrause}
R.~Hegselmann and U.~Krause.
\newblock Opinion dynamics and bounded confidence models, analysis, and
  simulation.
\newblock {\em J. Artif. Soc. Soc. Simul.}, 5(3):1--33, 2002.

\bibitem{HolmeSaramaki}
P.~Holme and J.~Saram{\"a}ki.
\newblock Temporal networks.
\newblock {\em Phys. Rep.}, 519(3):97--125, 2011.

\bibitem{KissMillerSimon}
I.Z. Kiss, J.~Miller, and P.L. Simon.
\newblock {\em Mathematics of Epidemics on Networks: From Exact to Approximate
  Models}.
\newblock Springer, 2017.

\bibitem{Kuramoto}
Y.~Kuramoto.
\newblock {\em Chemical Oscillations, Waves, and Turbulence}.
\newblock Springer, 2012.

\bibitem{Leicht2008}
E.~A. Leicht and M.~E.~J. Newman.
\newblock {Community Structure in Directed Networks}.
\newblock {\em Physical Review Letters}, 100(11):118703, mar 2008.

\bibitem{Newman}
M.E.J. Newman.
\newblock The structure and function of complex networks.
\newblock {\em SIAM Review}, 45:167--256, 2003.

\bibitem{PikovskyRosenblumKurths}
A.S. Pikovsky, M.~Rosenblum, and J.~Kurths.
\newblock {\em Synchronization}.
\newblock CUP, 2001.

\bibitem{Pons2005}
P.~Pons and M.~Latapy.
\newblock {Computing Communities in Large Networks Using Random Walks}.
\newblock In PInar Yolum, Tunga G{\"{u}}ng{\"{o}}r, Fikret G{\"{u}}rgen, and
  Can {\"{O}}zturan, editors, {\em Computer and Information Sciences - ISCIS
  2005}, pages 284--293, Berlin, Heidelberg, 2005. Springer Berlin Heidelberg.

\bibitem{PorterGleeson}
M.A. Porter and J.P. Gleeson.
\newblock {\em Dynamical Systems on Networks: A Tutorial}.
\newblock Frontiers in Applied Dynamical Systems: Reviews and Tutorials.
  Springer, 2016.

\bibitem{Reichardt2004}
J.~Reichardt and S.~Bornholdt.
\newblock {Detecting fuzzy community structures in complex networks with a
  potts model}.
\newblock {\em Physical Review Letters}, 93(21):19--22, 2004.

\bibitem{Reichardt2006}
J.~Reichardt and S.~Bornholdt.
\newblock {Statistical mechanics of community detection}.
\newblock {\em Physical Review E}, 74(1):016110, jul 2006.

\bibitem{Rossetti2019}
G.~Rossetti, L.~Milli, and R.~Cazabet.
\newblock {CDLIB: a python library to extract, compare and evaluate communities
  from complex networks}.
\newblock {\em Applied Network Science}, 4(1), 2019.

\bibitem{NetworkRepository}
R.~A. Rossi and N.~K. Ahmed.
\newblock The network data repository with interactive graph analytics and
  visualization.
\newblock In {\em AAAI}, 2015.

\bibitem{Traag2013}
V.~A. Traag, G.~Krings, and P.~{Van Dooren}.
\newblock {Significant Scales in Community Structure}.
\newblock {\em Scientific Reports}, 3(1):2930, dec 2013.

\bibitem{VicsekZafiris}
T.~Vicsek and A.~Zafiris.
\newblock Collective motion.
\newblock {\em Phys. Rep.}, pages 71--140, 2012.

\bibitem{WattsStrogatz}
D.J. Watts and S.H. Strogatz.
\newblock Collective dynamics of 'small-world' networks.
\newblock {\em Nature}, 393:440--442, 1998.

\end{thebibliography}

\newpage

\appendix

\section{Network Models}\label{sec:NetworkModels}

In this section we describe the pre-simulation steps (P1)-(P2) and the CIA Evaluation steps (E1)-(E2) in technical detail for several different network models. In a general network model of the form
\begin{align}\label{eq:GeneralNetwork}
    x_m'(t) = f_m(x_m) + \frac{1}{N} \sum_{\ell = 1}^N a_{m\ell}~g(x_\ell,x_m), \quad m=1,\dots,N
\end{align}
the steps (P2) and (E1) are quite abstract and general. When considering a specific network model, the function $g$ is often not just an abstract function but one that exhibits more structure that can be exploited in those two steps. The steps (P1) and (E2), however, do not depend on the specific network model. Even though we have only described (E2) in Section \ref{sec:CIA} for a specific example this step is exactly the same for all other network models.
Therefore, we shortly describe the community detection process (P1) in Subsection \ref{sec:CommunityDetection} such that we have covered the steps (P1) and (E2), which are mostly independent of the particular network model. In the following subsections we then do not further touch upon these steps but only describe (P2) and (E1) for specific network models. In particular, our fast CIAs work on block matrices with dense or fully occupied blocks. We shall explain below, why for illustration purposes, we may assume to illustrate the computation in the easiest case of just one full block. Hence, we assume for the subsections following Subsection \ref{sec:CommunityDetection} that $A$ is a full matrix with $A = (a_{m\ell})_{m,\ell=1,\dots,N}$ and $a_{m\ell} = 1$ for all $m,\ell$. Table \ref{tab:NetworkExamples} lists several examples of coupling functions of a wide variety of network models from different applications that are of the form \eqref{eq:GeneralNetwork}; in fact, there are many more application examples having a form identical or very similar to~\eqref{eq:GeneralNetwork} such as continuous Hopfield-type neural network models or the Hegselmann-Krause model for opinion formation. 

\begin{table}[h!]
\begin{center}
\begin{tabular}{|l|l|l|}\hline 
Network & State & Coupling Functions\\
Model & Space ($\mathcal X$) & \\ \hline

(DZ) & $\R$ & $f_m(x) = f(x) = - \, V'(x)$ \\
&& $h(\xi) = \xi$ \\ \hline 

(K) & $\S = \R/(2 \pi \Z)$ & $f_m(x) = \omega_m$ \\
& & $h(\xi) = \sin(\xi)$ \\ \hline 

(CS) & $\R^{2n}$ & $f_m(x) = (v,0)^\top$, $x = (s,v)$ \\
& & $h(\xi) = K~(\alpha^2 + \lVert \hat s\rVert_2^2)^{-\beta} \, (0, \hat v)^\top$, $\xi= (\hat s, \hat v)$ \\\hline

(V)/(FN) & $\R^2$
& $f_m(x) = f(x_1, x_2) = (x_2 - \tfrac{1}{3} \, x_1^3 + x_1, - \, \varepsilon \, x_1)^\top$ \\
& & $h(\xi) = \xi$ \\\hline

\end{tabular}

\caption{Desai--Zwanzig~(DZ), Kuramoto (K), Cucker--Smale~(CS), Van--der--Pol~(V) as well as FitzHugh--Nagumo~(FN) 
systems are relevant examples for time-continuous dynamical systems that can be cast into the form~\eqref{eq:GeneralNetworkSystem} with coupling functions $g(x,y) = h(x-y)$.}
\label{tab:NetworkExamples}
\end{center}
\end{table}

\subsection{Community Detection (P1)}\label{sec:CommunityDetection}

In general, our CIAs can quickly evaluate the right-hand side of \eqref{eq:GeneralNetwork}, when the interaction matrix $A = (a_{m\ell})_{m,\ell=1,\dots,N} \in \{0,1\}^{N\times N}$ is fully occupied by ones and so all nodes of the underlying graph form one large community. But CIAs can also quickly evaluate the right-hand side if the adjacency matrix has block structure such that the underlying graph is partitioned into smaller communities, within which each node is connected to every other node. Obviously, the adjacency matrices of most networks consist not only of blocks fully occupied by ones, but there may be a zero-entry of $A$ at a position $(m,\ell)$ that is in a larger block of ones. In other words, there is not always an all-to-all coupling inside communities but two nodes $m,\ell$ belonging to the same community might be uncoupled (missing intra-population links). Similarly, there may be a one-entry of $A$ at $(m,\ell)$ while nodes $m$ and $\ell$ do not belong the the same block in adjacency matrix $A$, meaning that two nodes $m$ and $\ell$ can be coupled even though they do not belong to the same community (additional inter-population links).
Our CIA nevertheless assumes that $A$ has block structure, calculates the right-hand side of \eqref{eq:GeneralNetwork} using this assumption in step (E1), and corrects this calculation by individually adding or subtracting $g(x_\ell, x_m)$ depending on whether $(m,\ell)$ is an additional inter-population link or a missing intra-population link of the graph in step (E2).

In synthetically created benchmark graphs, the communities of the graph, or the block structure of the associated adjacency matrix, respectively, is already known from the construction. However, when dealing with real-world networks, communities first have to be found.

CIAs can quickly evaluate the right-hand side for blocks in $A$ but the evaluation of many summands $g(x_\ell,x_m)$ is very costly. The aim of a community detection algorithm is partition the graph into communities such that there are as few missing intra-population links and additional inter-population links as possible. In other words, the matrix $S$ defined in Section \ref{sec:CIA} (see also Figure \ref{fig:Community_Detection}) has to be as sparse as possible. It can seen as a feature of a community detection algorithm to achieve exactly that. However, not all community detection algorithms pursue to optimize on that feature. Therefore, we have analyzed and compared different community detection algorithms with respect to that feature \cite{BohleKuehnThalhammer}.

Among the tested algorithms were \textsc{greedy\_modularity}~\cite{Clauset2004} from the \textsc{Python} package \textsc{networkx}~\cite{Hagberg2008} and the algorithms \textsc{louvain}~\cite{Blondel2008}, \textsc{rber\_pots}~\cite{Reichardt2004, Reichardt2006}, \textsc{rb\_pots}~\cite{Leicht2008, Reichardt2006}, \textsc{significance\_communities}~\cite{Traag2013}, \textsc{walktrap}~\cite{Pons2005} from the \textsc{Python} package \textsc{cdlib}~\cite{Rossetti2019}.

While all of these algorithms perform reasonably well, our tests showed that \textsc{rber\_pots} detects communities in a way such that the matrix $S$ has the fewest non-zero entries. Consequently this algorithm is best suitable for our applications. 
Knowing the theory behind this algorithm, it is no surprise that it performs well, since it is specifically optimized to reward existing intra-population links, punish missing intra-population links, reward missing inter-population links and punish existing inter-population links \cite{Reichardt2006}. This results in a Hamiltonian
\begin{align*}
    \mathcal H(\{\sigma\}) = -\sum_{m\neq \ell} (a_{m\ell} - \gamma p_{m\ell})\delta(\sigma_m, \sigma_\ell),
\end{align*}
that the algorithm tries to minimize \cite{Reichardt2006}. Here, $a_{m\ell}$ are the entries of the adjacency matrix of the graph and $p_{m\ell}$ is the probability that a link between node $m$ and $\ell$ exists. This probability is normalized such that $\sum{m\neq \ell} p_{m\ell} = 2N$. It can be chosen as $p_{m\ell} = k_mk_\ell / (2N)$, where $k$ denotes the degree distribution of the network, or one can just take $p_{m\ell} = |E|/((N-1)N)$, where $|E|$ is the total amount of edges in the network. Further, $\sigma_m$ is an index of the community to which node $m$ belongs such that $\delta(\sigma_m, \sigma_\ell) = 1$ if nodes $m$ and $\ell$ belong to the same community and $\delta(\sigma_m,\sigma_\ell) = 0$ otherwise. Finally, there is a parameter $\gamma$ that determines the ratio of how much a missing intra-population link should be punished in comparison with the reward of a non-existing inter-population link. In the standard case $\gamma = 1$, a missing intra-population link or an additional inter-population link negatively effects the Hamiltonian by the same amount as an existing intra-population link or a missing inter-population link positively effects it. Since we need to evaluate the coupling function $g(x_\ell, x_m)$ for each missing intra-population link and each existing inter-population link exactly once, $\gamma = 1$ is reasonable for our application.

Having praised \textsc{rber\_pots}, we also want to remark that all tested community detection algorithms do not take into account that the numerical evaluation of the right-hand side of \eqref{eq:FastHdense} takes some time, too. Since this evaluation time scales with the number of communities it is sometimes better to have fewer but larger communities, especially if the order $p$ of the expansion is high. 
The optimal community structure thus depends on the specific network model including its parameters such as the order $p$ of the expansion and the dimension of the model.


\subsection{General Network Model - Fourier Expansion}\label{sec:GeneralNetworkFourier}
Having established a community structure, here, we focus on just one community. For the sake of a clear notation, we take $(a_{m\ell})_{m,\ell=1,\dots,N} = (1)_{m,\ell=1,\dots,N}$, i.e. we consider the network model
\begin{align}\label{eq:GeneralNetworkFourier}
    x_m'(t) = \frac{1}{N} \sum_{\ell = 1}^N g(x_\ell,x_m), \quad m=1,\dots,N.
\end{align}
A naive computation of the right-hand side of \eqref{eq:GeneralNetworkFourier} for all $m=1,\dots,N$ requires $N^2$ evaluations of $g$ and thus the required time scales quadratically with $N$. Our aim is to reduce that to linear complexity.
For simplicity we first assume that the coupling function $g$ takes two inputs $x_\ell$ and $x_m$ from a one-dimensional space, e.g. the real line $\R$ or the circle $\S$ and maps into $\R$.

\paragraph{Complex Fourier Series} A general Fourier expansion of $g$ is then of the form
\begin{align}
    g(x,y) &= \sum_{\alpha \in \Z} \sum_{\beta \in \Z} c_{\alpha,\beta} e^{\frac{\pi i}{L}\alpha x}e^{\frac{\pi i}{L}\beta y}\\
    \label{eq:FourierApproximation}
    &\approx \sum_{\alpha = -p}^p \sum_{\beta =-p}^p c_{\alpha,\beta} e^{\frac{\pi i}{L}\alpha x}e^{\frac{\pi i}{L}\beta y},
\end{align}
where $c_{\alpha,\beta}$ are the Fourier coefficients of $g$, $L>0$ sets the size of the area $[-L,L]\times[-L,L]$ on which the Fourier expansion is valid and $p\in \N$ is an indicator for the approximation order. Inserting this representation into \eqref{eq:GeneralNetworkFourier}, we get
\begin{align*}
    x_m'(t) &= \frac{1}{N} \sum_{\ell = 1}^N g(x_\ell,x_m)\\
    &\approx \frac{1}{N} \sum_{\ell=1}^N \sum_{\alpha = -p}^p \sum_{\beta = -p}^p c_{\alpha,\beta}e^{\frac{\pi i}{L}\alpha x_\ell} e^{\frac{\pi i}{L}\beta x_m}\\
    &=\sum_{\alpha = -p}^p \sum_{\beta = -p}^p c_{\alpha,\beta}e^{\frac{\pi i}{L}\beta x_m}\underbrace{\left( \frac{1}{N} \sum_{\ell = 1}^N e^{\frac{\pi i}{L}\alpha x_\ell}\right)}_{=:r_\alpha}.
\end{align*}
Therefore, if one precomputes certain well-chosen observables, namely the so-called generalized order parameters 
\begin{align}\label{eq:FourierOrderParam}
    r_\alpha = \frac{1}{N} \sum_{\ell = 1}^N e^{\frac{\pi i}{L}\alpha x_\ell},
\end{align}
for all $\alpha = -p,\dots,p$, the computation of the right-hand side reduces to
\begin{align}\label{eq:FourierRHS}
    x_m'(t) \approx \sum_{\alpha = -p}^p \sum_{\beta = -p}^p c_{\alpha,\beta}\ r_\alpha \ e^{\frac{\pi i}{L}\beta x_m}.
\end{align}
Now, note that the precomputation complexity of the generalized order parameters \eqref{eq:FourierOrderParam} scales linearly in $N$ and so does the computation of \eqref{eq:FourierRHS}, since it has to be computed for all $m=1,\dots,N$. The total complexity thus has come down to $\mathcal O(N)$. The constant in front of the $N$ can be further improved by noting that $c_{\alpha,\beta} = \overline{c_{-\alpha,-\beta}}$ in \eqref{eq:FourierApproximation}, where $\bar c$ denotes the complex conjugate of $c$. This has to hold to guarantee that $g$ is a real function. Similarly, $r_\alpha = \overline{r_{-\alpha}}$ in \eqref{eq:FourierOrderParam}. Using these relations, the computation time can be further reduced by a constant factor, even if it still scales linearly with $N$. However, since these improvements would make the notation more elaborate and thus worsen the readability, we do not mention them further below. 

There still remains the question about how to choose $L$ and $p$ in the Fourier approximation. Unfortunately, there is no general answer to this question, but in specific cases the range of $x$ is restricted to a region $[-L,L]$ anyway, so then $L$ can be chosen such that the Fourier approximation is valid on the whole domain, see for example Section \ref{sec:Kuramoto} and \ref{sec:HOKuramoto}. Furthermore, in some cases, the coupling function is already a finite Fourier series, so \eqref{eq:FourierApproximation} is exact for small $p\in \N$ from which a choice of $p$ can be made. In general, $L$ should be chosen large enough such that $x_m(t) \in [-L,L]$ for all $k$ and all $t$ in the simulation time range. This can either numerically be tested or ensured by theoretical results that guarantee the boundedness of $x_m(t)$.

To summarize, the important steps are as follows:
\begin{itemize}
    \item[(P1)] Before starting the simulation, determine Fourier-coefficients $c_{\alpha,\beta}$, such that the coupling function $g$ can be represented as or well approximated by a finite series 
    \begin{align*}
        g(x,y) \approx \sum_{\alpha = -p}^p \sum_{\beta =-p}^p c_{\alpha,\beta} e^{\frac{\pi i}{L}\alpha x}e^{\frac{\pi i}{L}\beta y}.
    \end{align*}
    \item[(E1)] In each time step, precompute generalized order parameters
    \begin{align*}
        r_\alpha = \frac{1}{N} \sum_{\ell = 1}^N e^{\frac{\pi i}{L}\alpha x_\ell},
    \end{align*}
    for $\alpha = -p,\dots,p$ and calculate the right-hand side of \eqref{eq:GeneralNetworkFourier} based on the formula
    \begin{align*}
        x_m'(t) \approx \sum_{\alpha = -p}^p \sum_{\beta = -p}^p c_{\alpha,\beta}\ r_\alpha \ e^{\frac{\pi i}{L}\beta x_m}.
    \end{align*}
\end{itemize}

\paragraph{Real Fourier Series} Alternatively from the approach using complex approximations, we can also start with an approximation involving $\sin$ and $\cos$. Then, we first have to determine Fourier coefficients $c^{11}_{\alpha,\beta},c^{12}_{\alpha,\beta},c^{21}_{\alpha,\beta},c^{22}_{\alpha,\beta}$ such that

\begin{subequations}\label{eq:FourierRealExpansion}
\begin{align}
    g(x,y) \approx \sum_{\alpha = 0}^p \sum_{\beta = 0}^p \Bigg( &c^{11}_{\alpha,\beta} \cos\left(\frac{\pi}{L}\alpha x\right)\cos\left(\frac{\pi}{L}\beta y\right) + c^{12}_{\alpha,\beta}\cos\left(\frac{\pi}{L}\alpha x\right)\sin\left(\frac{\pi}{L}\beta y\right) \\&+ c^{21}_{\alpha,\beta} \sin\left(\frac{\pi}{L}\alpha x\right)\cos\left(\frac{\pi}{L}\beta y\right) + c^{22}_{\alpha,\beta} \sin\left(\frac{\pi}{L}\alpha x\right)\sin\left(\frac{\pi}{L}\beta y\right) \Bigg).
\end{align}
\end{subequations}

Using this representation, the right-hand side of \eqref{eq:GeneralNetworkFourier} reads as
\begin{align*}
    x_m'(t) &= \frac{1}{N} \sum_{\ell = 1}^N g(x_\ell,x_m) \\
    &\approx \frac{1}{N} \sum_{\ell = 1}^N \sum_{\alpha = 0}^p \sum_{\beta = 0}^p \Bigg( c^{11}_{\alpha,\beta} \cos\left(\frac{\pi}{L}\alpha x_\ell\right)\cos\left(\frac{\pi}{L}\beta x_m\right) + c^{12}_{\alpha,\beta}\cos\left(\frac{\pi}{L}\alpha x_\ell\right)\sin\left(\frac{\pi}{L}\beta x_m\right) \\&
    \qquad \qquad + c^{21}_{\alpha,\beta} \sin\left(\frac{\pi}{L}\alpha x_\ell\right)\cos\left(\frac{\pi}{L}\beta x_m\right) + c^{22}_{\alpha,\beta} \sin\left(\frac{\pi}{L}\alpha x_\ell\right)\sin\left(\frac{\pi}{L}\beta x_m\right) \Bigg)\\
    &=\sum_{\alpha = 0}^p \sum_{\beta = 0}^p \Bigg[ c^{11}_{\alpha,\beta} \left( \frac 1N \sum_{\ell=1}^N \cos\left(\frac{\pi}{L}\alpha x_\ell\right)\right) \cos\left(\frac{\pi}{L}\beta x_m\right) + c^{12}_{\alpha,\beta} \left( \frac 1N \sum_{\ell=1}^N \cos\left(\frac{\pi}{L}\alpha x_\ell\right)\right) \sin\left(\frac{\pi}{L}\beta x_m\right) \\& \qquad \qquad  +c^{21}_{\alpha,\beta} \left( \frac 1N \sum_{\ell=1}^N \sin\left(\frac{\pi}{L}\alpha x_\ell\right)\right) \cos\left(\frac{\pi}{L}\beta x_m\right) + c^{22}_{\alpha,\beta} \left(\frac 1N \sum_{\ell=1}^N \sin\left(\frac{\pi}{L}\alpha x_\ell\right)\right)\sin\left(\frac{\pi}{L}\beta x_m\right) \Bigg].
\end{align*}
Therefore, in each time step, we need to precompute 
\begin{align}\label{eq:RealOP}
    r_\alpha^\text{cos} := \frac 1N \sum_{\ell=1}^N \cos\left(\frac{\pi}{L}\alpha x_\ell\right) \quad \text{and} \quad r_\alpha^\text{sin} := \frac 1N \sum_{\ell=1}^N \sin\left(\frac{\pi}{L}\alpha x_\ell\right)
\end{align}
for all $\alpha = 0,\dots,p$. Having done that, the right-hand side of \eqref{eq:GeneralNetworkFourier} can be rewritten as
\begin{align*}
    x_m'(t) &\approx \sum_{\alpha = 0}^p \sum_{\beta = 0}^p \Bigg[ c^{11}_{\alpha,\beta} r_\alpha^\text{cos} \cos\left(\frac{\pi}{L}\beta x_m\right) + c^{12}_{\alpha,\beta} r_\alpha^\text{cos} \sin\left(\frac{\pi}{L}\beta x_m\right) \\& \qquad \qquad  +c^{21}_{\alpha,\beta} r_\alpha^\text{sin} \cos\left(\frac{\pi}{L}\beta x_m\right) + c^{22}_{\alpha,\beta} r_\alpha^\text{sin}\sin\left(\frac{\pi}{L}\beta x_m\right) \Bigg].
\end{align*}

To summarize, the important steps when using a real expansion are
\begin{itemize}
    \item[(P1)] Before starting the simulation, determine Fourier-coefficients $c^{11}_{\alpha,\beta},c^{12}_{\alpha,\beta},c^{21}_{\alpha,\beta},c^{22}_{\alpha,\beta}$, such that the coupling function $g$ can be represented as or well approximated by a finite series
        \begin{align}
            g(x,y) \approx \sum_{\alpha = 0}^p \sum_{\beta = 0}^p \Bigg( &c^{11}_{\alpha,\beta} \cos\left(\frac{\pi}{L}\alpha x\right)\cos\left(\frac{\pi}{L}\beta y\right) + c^{12}_{\alpha,\beta}\cos\left(\frac{\pi}{L}\alpha x\right)\sin\left(\frac{\pi}{L}\beta y\right) \\&+ c^{21}_{\alpha,\beta} \sin\left(\frac{\pi}{L}\alpha x\right)\cos\left(\frac{\pi}{L}\beta y\right) + c^{22}_{\alpha,\beta} \sin\left(\frac{\pi}{L}\alpha x\right)\sin\left(\frac{\pi}{L}\beta y\right) \Bigg).
        \end{align}
    \item[(E1)] In each time step, precompute
        \begin{align*}
            r_\alpha^\text{cos} := \frac 1N \sum_{\ell=1}^N \cos\left(\frac{\pi}{L}\alpha x_\ell\right) \quad \text{and} \quad r_\alpha^\text{sin} := \frac 1N \sum_{\ell=1}^N \sin\left(\frac{\pi}{L}\alpha x_\ell\right)
        \end{align*}
        for $\alpha,\beta = 0,\dots,p$ and calculate the right-hand side of \eqref{eq:GeneralNetworkFourier} based on the formula
        \begin{align*}
            x_m'(t) &\approx \sum_{\alpha = 0}^p \sum_{\beta = 0}^p \Bigg[ c^{11}_{\alpha,\beta} r_\alpha^\text{cos} \cos\left(\frac{\pi}{L}\beta x_m\right) + c^{12}_{\alpha,\beta} r_\alpha^\text{cos} \sin\left(\frac{\pi}{L}\beta x_m\right) \\& \qquad \qquad  +c^{21}_{\alpha,\beta} r_\alpha^\text{sin} \cos\left(\frac{\pi}{L}\beta x_m\right) + c^{22}_{\alpha,\beta} r_\alpha^\text{sin}\sin\left(\frac{\pi}{L}\beta x_m\right) \Bigg].
        \end{align*}
\end{itemize}

\paragraph{Difference based coupling function - Complex Fourier Series} Even though we have already reduced the complexity from being quadratic in $N$ to being only linear in $N$, the constant scales with $p^2$. In many network models the coupling function $g$ is of the special form $g(x,y) = h (x-y)$, such that we are facing the system
\begin{align}\label{eq:GeneralNetworkFourierDifference}
    x_m'(t) = \frac{1}{N} \sum_{\ell = 1}^N h(x_\ell-x_m), \quad k=1,\dots,N.
\end{align}
This helps to reduce the dependence on $p^2$ to just $p$. Again, we assume that the coupling function $h$ has an Fourier approximation
\begin{align}\label{eq:FourierDifferenceApproximation}
    h(x) &= \sum_{\alpha\in \Z} d_\alpha e^{\frac{\pi}{L}\alpha x} \approx \sum_{\alpha = -p}^p d_\alpha e^{\frac{\pi}{L}\alpha x},
\end{align}
where $d_\alpha$ are the Fourier coefficients, $L>0$ indicates the size of the domain $[-L,L]$ in which the approximation is valid and $p$ gives the approximation order. Then, the right-hand side of \eqref{eq:GeneralNetworkFourier} can be written as
\begin{align}
    \nonumber
    x_m'(t) & = \frac{1}{N} \sum_{\ell=1}^N h(x_\ell-x_m)\\
    \nonumber
    &\approx \frac 1N \sum_{\ell=1}^N \sum_{\alpha = -p}^p d_\alpha e^{\frac{\pi}{L} \alpha(x_\ell-x_m)}\\
    &=\sum_{\alpha = -p}^pd_\alpha \left( \frac 1N \sum_{\ell=1}^N e^{\frac \pi L \alpha x_\ell}\right) e^{-\frac \pi L \alpha x_m}\\
    \label{eq:FourierRHSDifferenceFast}
    &=\sum_{\alpha = -p}^pd_\alpha r_\alpha e^{-\frac \pi L \alpha x_m},
\end{align}
where $r_\alpha$ are the generalized order parameters \eqref{eq:FourierOrderParam}. While the general formula \eqref{eq:FourierRHS} involves two sums with indices running from $-p$ to $p$, the formula \eqref{eq:FourierRHSDifferenceFast}, which relies on the assumption of a difference based coupling, involves only one such sum.

\paragraph{Difference based coupling function - Real Fourier Series}
Again, instead of expanding $h$ in a complex Fourier series, one can also use a real Fourier series
\begin{align}\nonumber
    h(x) &= d_0^\text{cos} + \sum_{\alpha = 1}^\infty \left( d_\alpha^\text{sin} \sin\left( \frac \pi L \alpha x\right) + d_\alpha^\text{cos} \cos\left( \frac \pi L \alpha x \right) \right)\\
    \label{eq:FourierRealDifferenceApproximation}
    &\approx d_0^\text{cos} + \sum_{\alpha = 1}^p \left( d_\alpha^\text{sin} \sin\left( \frac \pi L \alpha x\right) + d_\alpha^\text{cos} \cos\left( \frac \pi L \alpha x \right) \right)
\end{align}
Having precomputed the quantities $r_\alpha^\text{sin}$ and $r_\alpha^\text{cos}$ from \eqref{eq:RealOP}, the right-hand side of \eqref{eq:GeneralNetworkFourierDifference} is given by
\begin{align*}
    x_m'(t) &= \frac{1}{N} \sum_{\ell = 1}^N h(x_\ell-x_m)\\
    &\approx \frac{1}{N} \sum_{\ell = 1}^N \left[d_0^\text{cos} + \sum_{\alpha = 1}^p \left( d_\alpha^\text{sin} \sin\left( \frac \pi L \alpha (x_\ell-x_m)\right) + d_\alpha^\text{cos} \cos\left( \frac \pi L \alpha (x_\ell-x_m) \right) \right)\right]\\
    &=d_0^\text{cos} + \frac 1N \sum_{\ell=1}^N \sum_{\alpha = 1}^p\Bigg[ d^\text{sin}_\alpha \sin\left(\frac \pi L \alpha x_\ell\right)\cos\left(\frac \pi L \alpha x_m\right) - d^\text{sin}_\alpha \cos\left(\frac \pi L \alpha x_\ell \right)\sin(\left(\frac \pi L \alpha x_m \right)\\
    &\qquad \qquad +d^\text{cos}_\alpha \sin\left( \frac \pi L \alpha x_\ell\right) \sin\left(\frac \pi L \alpha x_m\right) + d^\text{cos}_\alpha \cos\left( \frac \pi L \alpha x_\ell \right) \cos\left( \frac \pi L \alpha x_m\right) \Bigg]\\
    &= d_0^\text{cos}+ \sum_{\alpha = 1}^p\Bigg[ d^\text{sin}_\alpha r^\text{sin}_\alpha \cos\left(\frac \pi L \alpha x_m\right) - d^\text{sin}_\alpha r_\alpha^\text{cos} \sin\left(\frac \pi L \alpha x_m\right)\\
    &\qquad \qquad + d^\text{cos}_\alpha r_\alpha^\text{sin} \sin\left(\frac \pi L \alpha x_m\right) + d^\text{cos}_\alpha r_\alpha^\text{cos} \cos\left( \frac \pi L \alpha x_m\right) \Bigg].
\end{align*}
This last equation represents a the formula that one should use to compute the right-hand side of \eqref{eq:GeneralNetworkFourierDifference} when preferring real Fourier approximations.

\paragraph{Extensions} In the above calculations we assumed that $g$ or $h$ take inputs form a one-dimensional space and map into a one-dimensional space. However, we want to remark that this approach also works if the inputs $x_\ell$ and $x_m$ are higher-dimensional objects, for example, when $g\colon \R^n\times \R^n \to \R^n$. In this case, $\alpha$ and $\beta$ have to be thought of being multi-indices rather than integers. Quantities of the form $e^{\frac{\pi i}{L} \alpha x}$ have to be replaced with $e^{\frac {\pi i}{L} \langle \alpha, x\rangle}$, where $\langle\cdot,\cdot\rangle$ is a scalar product. Further, sums over $\alpha,\beta = -p,\dots,p$ are then sums over $\alpha,\beta \in \Z(p)^n:=\{-p,\dots,p\}^n$. Moreover, the order parameter \eqref{eq:FourierOrderParam} or its real equivalents need to be precomputed for all $\alpha\in \Z(p)^n$.

\subsection{General Network Model - Polynomial Expansion}\label{sec:GeneralNetworkPolynom}
Again, we consider the general network model
\begin{align}\label{eq:GeneralNetworkPoly}
    x_m'(t) = \frac{1}{N} \sum_{\ell = 1}^N g(x_\ell,x_m), \quad m=1,\dots,N.
\end{align}
Again, our goal is to reduce the computational complexity from $N^2$ to just $N$. We assume for simplicity that $g$ takes two inputs from a one-dimensional space such as $\R$ or $\S$ and maps to $\R$. However, it should be said that our approach works as well when the inputs of $g$ are from a higher-dimensional space. 

\paragraph{Polynomial Approximation} However, instead of approximating the coupling function by a Fourier series, this time we approximate it by polynomials
\begin{align}\label{eq:PolyApprox}
    g(x,y) = \sum_{\alpha = 0}^\infty \sum_{\beta = 0}^\infty c_{\alpha,\beta} x^\alpha y^\beta \approx \sum_{\alpha = 0}^p \sum_{\beta = 0}^p c_{\alpha,\beta} x^\alpha y^\beta.
\end{align}
Here $c_{\alpha,\beta}$ are the coefficients of the approximation and $p\in \N$ indicates the approximation order. This approximation does not necessarily need to be a Taylor approximation. Rather, it is often more useful to consider a polynomial approximation of $g$ with respect to a $\mathrm L^2$ or a supremum norm on a domain $[-L,L] \times [-L,L]$. For numerical reasons it sometimes make sense to replace $x$ and $y$ in \eqref{eq:PolyApprox} by $(x-x_0)$ and $(y-y_0)$, respectively.
Combining this approximation with \eqref{eq:GeneralNetworkPoly}, we obtain
\begin{align*}
    x_m'(t) &= \frac{1}{N} \sum_{\ell = 1}^N g(x_\ell,x_m)\\
    &\approx \frac 1N \sum_{\ell=1}^N \sum_{\alpha = 0}^p\sum_{\beta = 0}^p c_{\alpha,\beta} x_\ell^\alpha x_m^\beta\\
    &= \sum_{\alpha = 0}^p\sum_{\beta = 0}^p c_{\alpha,\beta} \underbrace{\left( \frac 1N \sum_{\ell=1}^N x_\ell^\alpha \right)}_{=:w_\alpha} x_m^\beta.
\end{align*}
Therefore, if one precomputes the $\alpha$-th moments
\begin{align}\label{eq:PolyMoments}
    w_\alpha := \frac 1N \sum_{\ell=1}^N x_\ell^\alpha
\end{align}
for $\alpha = 0,\dots,p$, the computation of the right-hand side reduces to
\begin{align}\label{eq:PolyRHS}
    x_m'(t) \approx \sum_{\alpha = 0}^p\sum_{\beta = 0}^p c_{\alpha,\beta} w_\alpha x_m^\beta.
\end{align}

To summarize, the important steps are as follows:
\begin{itemize}
    \item[(P1)] Before starting the simulation, determine coefficients $c_{\alpha,\beta}$, such that the coupling function $g$ can be represented as or well approximated by a finite series
    \begin{align*}
        g(x,y) \approx \sum_{\alpha = 0}^p \sum_{\beta = 0}^p c_{\alpha,\beta} x^\alpha y^\beta.
    \end{align*}
    \item[(E1)] In each time step, precompute the moments
    \begin{align*}
         w_\alpha := \frac 1N \sum_{\ell=1}^N x_\ell^\alpha
    \end{align*}
    for $\alpha = 0,\dots,p$ and calculate the right-hand side of \eqref{eq:GeneralNetworkPoly} based on the formula
    \begin{align*}
        x_m'(t) \approx \sum_{\alpha = 0}^p\sum_{\beta = 0}^p c_{\alpha,\beta} w_\alpha x_m^\beta.
    \end{align*}
\end{itemize}
As one can see, the complexity of an evaluation of the right-hand side \eqref{eq:PolyRHS} is only linear in $N$, since it has to be evaluated for each $k=1,\dots,N$. The dependence of this complexity on $p^2$ can be reduced in special cases, for example if the coupling function $g$ depends only on differences.

\paragraph{Difference based coupling function - Polynomial Approximation}
Even though we have already reduces the complexity from being quadratic in $N$ to being only linear in $N$, the constant scales with $p^2$. In many network models the coupling function $g$ is of the special form $g(x,y) = h (x-y)$, such that we are facing the system
\begin{align}\label{eq:GeneralNetworkPolyDifference}
    x_m'(t) = \frac{1}{N} \sum_{\ell = 1}^N h(x_\ell-x_m), \quad k=1,\dots,N.
\end{align}
Such a representation is helpful when one wants to further reduce the computational complexity. Now, we assume that the coupling function $h$ can be well approximated by a polynomial
\begin{align*}
    h(x) = \sum_{\alpha=0}^\infty c_\alpha x^\alpha \approx \sum_{\alpha=0}^p c_\alpha x^\alpha.
\end{align*}
Again, for numerical reasons it is sometimes better to replace $x$ in the above formula with $(x-x_0)$. However, for the sake of simplicity we do not incorporate this technical detail. Given this polynomial approximation and the $\alpha$-th moments \eqref{eq:PolyMoments}, we can rewrite the right-hand side of \eqref{eq:GeneralNetworkPolyDifference} to
\begin{align}
    \nonumber
    x_m'(t) &= \frac{1}{N} \sum_{\ell = 1}^N h(x_\ell-x_m)\\
    \nonumber
    &\approx \frac{1}{N} \sum_{\ell = 1}^N \sum_{\alpha = 0}^p c_\alpha (x_\ell-x_m)^\alpha\\
    \nonumber
    &= \sum_{\alpha = 0}^p c_\alpha \frac 1N \sum_{\ell=1}^N \sum_{k = 0}^\alpha \begin{pmatrix} \alpha\\ k\end{pmatrix} x_\ell^k ~ (-x_m)^{\alpha-k}\\
    \nonumber
    &=\sum_{\alpha = 0}^p c_\alpha \sum_{k = 0}^\alpha \begin{pmatrix} \alpha\\ k\end{pmatrix} \left( \frac 1N \sum_{\ell=1}^N x_\ell^k\right) (-x_m)^{\alpha-k}\\
    \label{eq:PolyDifferenceRHSfast}
    &=\sum_{\alpha = 0}^p c_\alpha \sum_{k = 0}^\alpha \begin{pmatrix} \alpha\\ k\end{pmatrix} w_k ~(-x_m)^{\alpha-k}.
\end{align}
This representation further reduces the computational complexity.

\subsection{Kuramoto Model}\label{sec:Kuramoto}

The classical Kuramoto model~\cite{Kuramoto} is given by
\begin{align*}
    \frac \dd {\dd t} \theta_m(t) = \omega_m + \frac 1N \sum_{\ell=1}^N \sin(\theta_\ell(t)-\theta_m(t)),\quad m=1,\dots,N,
\end{align*}
where $\theta_m\colon [0,T]\to \S = \R / (2\pi \Z)$ is the phase and $\omega_m$ the intrinsic frequency of oscillator $m$. The coupling function $g$ is hence given by $g(x,y) = h (x-y) = \sin(x-y)$.

\paragraph{Classical Kuramoto model}
Following the difference based approach from Section \ref{sec:GeneralNetworkFourier}, we choose $L=\pi$, $p=1$ and write $h$ as
\begin{align*}
    h (x) = \sin(x) = \frac{-1}{2i}e^{-ix} + \frac{1}{2i} e^{ix},
\end{align*}
so $d_{-1} = -1/(2i), d_0 = 0$ and $d_1=1/(2i)$ in \eqref{eq:FourierDifferenceApproximation} and this approximation is exact.
After calculating the order parameters
\begin{align}
    \nonumber
    r_{-1} &= \frac 1N \sum_{\ell=1}^N e^{-i\theta_\ell},\\
    \nonumber
    r_0 &= 1,\\
    \label{eq:KuramotoFirstOP}
    r_1 &= \frac 1N \sum_{\ell=1}^N e^{i\theta_\ell},
\end{align}
the equation to evaluate the right-hand side \eqref{eq:FourierRHSDifferenceFast} turns into
\begin{align}
    \label{eq:Kuramotofast1}
    \theta_m'(t) &= \omega_m + \frac{-1}{2i} r_{-1} e^{i\theta_m} + \frac{1}{2i} r_1 e^{-i\theta_m}.
\end{align}
However, since $d_{-1} = \overline{d_1}$ and $r_{-1} = \overline{r_1}$ we can further simplify:
\begin{align}
    \label{eq:Kuramotofast2}
    \theta_m'(t) &= \omega_m - \operatorname{Re}\left(i r_1 e^{-i\theta_m}\right)\\
    \label{eq:Kuramotofast3}
    &= \omega_m + \operatorname{Im}(r_1 e^{-i\theta_m}).
\end{align}
One can also write $r_1 = \abs{r_1} e^{i\psi}$ for some $\psi\in \S$. Then,
\begin{align}
    \nonumber
    \theta_m'(t) &= \omega_m + \operatorname{Im}(r_1 e^{-i \theta_m})\\
    \label{eq:Kuramotofast4}
    &=\omega_m + \abs{r_1} \operatorname{Im}(e^{i (\psi-\theta_m)})\\
    \label{eq:Kuramotofast5}
    &=\omega_m + \abs{r_1}\sin(\psi - \theta_m).
\end{align}

Alternatively, one can also prefer to work with real numbers only. Then, one has to precompute
\begin{align}
    \label{eq:KuramotoFirstOPReal}
    r^\text{cos}_1 = \frac{1}{N}\sum_{\ell=1}^N \cos(\theta_\ell),\quad  \text{and} \quad r^\text{sin}_1 = \frac 1N \sum_{\ell=1}^N \sin(\theta_\ell).
\end{align}
According to the derivation in Section \ref{sec:GeneralNetworkFourier}, for the right-hand side we obtain
\begin{align}
    \label{eq:Kuramotofast6}
    \theta_m'(t) &= \omega_m + r^\text{sin}_1\cos(x_m) - r^\text{cos}_1\sin(x_m).
\end{align}

To summarize, the important steps are given by
\begin{itemize}
    \item[(P1)] In this step nothing has to be done, since the coupling function is already a finite Fourier series.
    \item[(E1)] In each simulation step, first calculate the complex order parameter \eqref{eq:KuramotoFirstOP} and then evaluate the right-hand side by using either of the formulas \eqref{eq:Kuramotofast1},\eqref{eq:Kuramotofast2},\eqref{eq:Kuramotofast3},\eqref{eq:Kuramotofast4},\eqref{eq:Kuramotofast5}. Alternatively, calculate the real order-parameters \eqref{eq:KuramotoFirstOPReal} and then evaluate the right-hand side using the formula \eqref{eq:Kuramotofast6}.
\end{itemize}

\paragraph{Higher-Harmonics Kuramoto Model}
An easy generalization of the classical Kuramoto model additionally includes higher harmonics in the coupling function. The network model is then given by
\begin{align*}
    \theta_m'(t) = \omega_m + \frac 1N \sum_{\ell=1}^N h(\theta_\ell(t)-\theta_m(t)),
\end{align*}
with a coupling function $h\colon \S \to \R$ defined by
\begin{align*}
    h(x) = \sum_{\alpha = 1}^p ( d_\alpha^\text{sin} \sin(\alpha x) + d_\alpha^\text{cos} \cos(\alpha x) ).
\end{align*}
Obviously, this is already of the form \eqref{eq:FourierRealDifferenceApproximation} for $L=\pi$, so we can directly follow this section. Having precomputed the quantities $r^\text{cos}_\alpha$ and $r^\text{sin}_\alpha$ from  \eqref{eq:RealOP} for all $\alpha = 0,\dots,p$, the right-hand side for the Higher-Harmonics Kuramoto Model can be written as
\begin{align*}
    \theta_m'(t) &= \omega_k + d_0^\text{cos} + \sum_{\alpha = 1}^p\Big[ d^\text{sin}_\alpha r^\text{sin}_\alpha \cos\left(\alpha x_m\right) - d^\text{sin}_\alpha r_\alpha^\text{cos} \sin\left( \alpha x_m\right)\\
    &\qquad \qquad\qquad  + d^\text{cos}_\alpha r_\alpha^\text{sin} \sin\left( \alpha x_m\right) + d^\text{cos}_\alpha r_\alpha^\text{cos} \cos\left(  \alpha x_m\right) \Big].
\end{align*}

\subsection{Desai-Zwanzig Model}\label{sec:DesaiZwanzig}

The Desai-Zwanzig model~\cite{DesaiZwanzig} is given by the following set of equations: 
\begin{align}\label{eq:DesaiZwanzig}
    x_m'(t) = -V'(x_m) + \frac{1}{N}\sum_{\ell=1}^N (x_\ell-x_m),
\end{align}
where $V\colon \R\to \R$ is a potential. Following the difference based procedures in Section \ref{sec:GeneralNetworkPolynom}, the coupling function $h$ is given by just $h(x) = x$. Therefore, to match the notation in this section, $p=1, c_0=0$ and $c_1=1$. After having computed the first moment $w_1$ from \eqref{eq:PolyMoments} and by using \eqref{eq:PolyDifferenceRHSfast}, we can write the right-hand side as
\begin{align*}
    x_m'(t) = -V'(x_m) + (-x_m + w_1).
\end{align*}
Even though this is an easy application of the theory from Section \ref{sec:GeneralNetworkPolynom} and could have easily derived from \eqref{eq:DesaiZwanzig} without the general theory from this section, it helps to reduce the computational cost significantly and thereby lowers the complexity from $\mathcal O(N^2)$ to just $\mathcal O(N)$.

\subsection{Cucker-Smale Model}\label{sec:CuckerSmale}

The continuous Cucker-Smale model~\cite{CuckerSmale1} is given by the dynamical system
\begin{subequations}
\label{eq:Cucker_Smale}
\begin{align}
    \label{eq:Cucker_Smale_Position}
    s_m'(t) &= v_m(t)\\
    \label{eq:Cucker_Smale_Velocity}
    v_m'(t) &= \frac{1}{N} \sum_{\ell=1}^N \frac{K}{(\sigma^2 + \lVert s_\ell(t)-s_m(t)\rVert^2)^\beta}~(v_\ell(t)-v_m(t)),
\end{align}
\end{subequations}
for $m=1,\dots,N$, where $s_m(t)$ represents the current position of the $m$-th bird, $v_m(t)$ is its velocity and $K, \sigma$ and $\beta$ are coupling constants. Here $s_m(t), v_m(t)\in \R^n$, where typically $n = 1,2,3$.
Putting this into the form \eqref{eq:GeneralNetwork} with $x = (s,v)^\top$, $f_m(x):=(v,0)^\top$ would be functions mapping from $\R^{2n}$ to $\R^{2n}$ and similarly, $g\colon \R^{2n}\times\R^{2n}\to \R^{2n}$, with $g(\hat x, \tilde x) = (g_s(\hat x, \tilde x), g_v(\hat x,\tilde x))^\top$. Here, the first $n$ components of $g$ are given by $g_s(\hat x, \tilde x) = 0$ and the last $n$ components of $g$ are given by 
\begin{align*}
    g_v\left(\hat x,\tilde x\right) = \eta(\lVert \hat s - \tilde s\rVert^2) ~ (\hat v - \tilde v),\quad \text{with } \eta(y) = \frac{K}{(\sigma^2+y)^\beta}.
\end{align*}
If we directly applied the algorithm described in previous section \ref{sec:GeneralNetworkFourier} or \ref{sec:GeneralNetworkPolynom}, the high dimension of the Cucker-Smale model would impact the performance of these algorithms, since they do not account for the special structure of the model. However, by exploiting this special structure, a more efficient algorithm can be constructed. In particular, a more efficient algorithm has to take into account that $g_s = 0$, so there is no need to expand this part either in a Fourier or a polynomial series. Furthermore, $g(\hat x, \tilde x)$ only depends on the difference $\hat x-\tilde x$, which should be exploited. Moreover, $g_v(\hat x, \tilde x)$ depends on $\hat v-\tilde v$ only linearly, so a polynomial expansion up to a degree higher than $1$ is unnecessary. Last but not least, the fraction in the sum of \eqref{eq:Cucker_Smale_Velocity} is independent of the coordinate direction, which makes it superfluous to expand this fraction for each coordinate direction.

Let us start developing a fast algorithm by denoting $\tilde \eta(y) \colon \R^n\to \R$, with $\tilde \eta(y) := \eta( \lVert y\rVert ^2)$ and expanding this in a Fourier series
\begin{align}
    \label{eq:Cucker_Smale_etaExpansion}
    \tilde \eta(y) = \sum_{\alpha\in \Z^n} \tilde c_\alpha e^{\frac{i\pi}L\langle\alpha, y\rangle} \approx \sum_{\alpha\in \Z(p)^n} \tilde c_\alpha e^{\frac{i\pi}L\langle\alpha, y\rangle},
\end{align}
where $\alpha\in \Z(p)^n\subset \Z^n$ is a multi-index, $\Z(p) = \{-p,\dots,p\}$, $\tilde c_\alpha$ are the Fourier coefficients of $\tilde \eta$, $L>0$ is a parameter that denotes the region $[-L,L]^n$ on which the expansion is valid and $\langle \alpha, x\rangle = \sum_{\ell = 1}^n \alpha_\ell x_\ell$ denotes the standard scalar product. In this new notation the second component of the right-hand side of \eqref{eq:Cucker_Smale_Velocity} reads as
\begin{align*}
    v_m'(t) = \frac{1}{N} \sum_{\ell=1}^N \tilde \eta (s_\ell(t)-s_m(t))~(v_\ell(t)-v_m(t)).
\end{align*}
Inserting the approximation \eqref{eq:Cucker_Smale_etaExpansion} into this formula yields
\begin{align*}
    v_m'(t) &\approx \frac{1}{N} \sum_{\ell=1}^N \sum_{\alpha\in \Z(p)^n} \tilde{c}_\alpha e^{\frac {i\pi}L \langle \alpha, s_\ell(t)-s_m(t)\rangle} (v_\ell(t)-v_m(t))\\
    &= \sum_{\alpha\in \Z(p)^n}\tilde{c}_\alpha e^{\frac {i\pi}{L} \langle \alpha, - s_m(t)\rangle} \left( \frac{1}{N} \sum_{\ell=1}^N e^{\frac {i\pi}L \langle \alpha, s_\ell(t)\rangle} v_\ell(t) - \frac{1}{N} \sum_{\ell=1}^N e^{\frac {i\pi}L \langle \alpha, s_\ell(t)\rangle}v_m(t)  \right).
\end{align*}
Therefore, if one precomputes
\begin{subequations}
\label{eq:Cucker_Smale_orderParams}
\begin{align}
    u_\alpha &= \frac 1N \sum_{\ell=1}^N e^{\frac {i\pi}L \langle \alpha, s_\ell\rangle} \in \R,\\
    h_\alpha &= \frac 1N \sum_{\ell=1}^N e^{\frac {i\pi}L \langle \alpha, s_\ell\rangle}v_\ell \in \R^n
\end{align}
\end{subequations}
for each $\alpha \in \Z(p)^n$, the right-hand side is finally given by
\begin{align}\label{eq:Cucker_Smale_fast_RHS}
    v_m'(t) \approx \sum_{\alpha\in \Z(p)^n} \tilde c_\alpha \ e^{\frac {i\pi}{L} \langle \alpha, -s_m(t)\rangle} \left( h_\alpha - u_\alpha v_m(t)\right).
\end{align}

To summarize, the important steps are given by
\begin{itemize}
    \item[(P1)] Before starting the simulation, determine Fourier-coefficients $\tilde c_{\alpha}$, such that the function $\tilde \eta$ is well approximated by a finite series of the form \eqref{eq:Cucker_Smale_etaExpansion}.
    \item[(E1)] In each time step, precompute the quantities \eqref{eq:Cucker_Smale_orderParams} for all $\alpha \in \Z(p)^n$ and calculate the right-hand side of \eqref{eq:Cucker_Smale_Velocity} by using the formula \eqref{eq:Cucker_Smale_fast_RHS}.
\end{itemize}

\subsection{Higher-order Kuramoto Models}\label{sec:HOKuramoto}

Higher-order Kuramoto models are generalizations from the classical Kuramoto model. While in the classical Kuramoto model the particle interactions are pairwise, in higher-order Kuramoto models, interactions of triplets, quadruplets, etc. determine the dynamics. We are going to show below that the CIA approach naturally generalizes to higher-order coupling, and is even more powerful in this case. The higher-order Kuramoto network model we use as an illustration for this generalization is given by
\begin{align*}
    \theta_m'(t) = \frac{1}{N^n}\sum_{\ell\in [N]^n} \sin\left( \sum_{k = 1}^n \lambda_k \theta_{\ell_k}(t) + \lambda_{n+1}\theta_m(t)\right).
\end{align*}
Here, $N$ is the number of oscillators, $[N]=\{1,\dots,N\}$,  $n + 1$ is amount of oscillators that interact with each other, $\theta_m(t)$ are the phases of oscillators $m=1,\dots,N$ and $\lambda_1,\dots,\lambda_{n+1}\in \Z$ are integer valued coefficients that typically sum up to $0$. In the classical Kuramoto model $n=1$, $\lambda_1 = 1, \lambda_2 = -1$. To make the presentation simple, we restrict ourself to the model higher-order model
\begin{align}\label{eq:HOKuramoto}
    \theta_m'(t) = \frac{1}{N^3}\sum_{\ell_1,\ell_2,\ell_3=1}^N \sin \Big(\lambda_1\theta_{\ell_1}(t) + \lambda_2\theta_{\ell_2}(t) + \lambda_3\theta_{\ell_3}(t) + \lambda_4\theta_m(t) \Big).
\end{align}
The computational complexity to naively evaluate the right-hand side of \eqref{eq:HOKuramoto} is $\mathcal O(N^4)$ since there are $3$ sums and they have to be evaluated for each $m=1,\dots,N$. The following calculation shows how to reduce this complexity:
\begin{align}
    \nonumber
    \theta_m'(t) &= \frac{1}{N^3}\sum_{\ell_1,\ell_2,\ell_3=1}^N \sin \Big(\lambda_1\theta_{\ell_1}(t) + \lambda_2\theta_{\ell_2}(t) + \lambda_3\theta_{\ell_3}(t) + \lambda_4\theta_m(t) \Big)\\
    \nonumber
    &=\operatorname{Im}\left( \frac{1}{N^3}\sum_{\ell_1,\ell_2,\ell_3=1}^N e^{i (\lambda_1\theta_{\ell_1}(t) + \lambda_2\theta_{\ell_2}(t) + \lambda_3\theta_{\ell_3}(t) + \lambda_4\theta_m(t))}\right)\\
    \nonumber
    &=\operatorname{Im}\left( \frac{1}{N^3} \sum_{\ell_1,\ell_2,\ell_3=1}^N e^{i \lambda_1\theta_{\ell_1}(t)}~ e^{i \lambda_2\theta_{\ell_2}(t)}~ e^{i \lambda_3\theta_{\ell_3}(t)}~ e^{i \lambda_4\theta_m(t)} \right)\\
    \nonumber
    &=\operatorname{Im}\left(\left(\frac 1N \sum_{j_1=1}^N e^{i \lambda_1\theta_{\ell_1}(t)}\right) \left(\frac 1N \sum_{j_2=1}^N e^{i \lambda_2\theta_{\ell_2}(t)}\right) \left(\frac 1N \sum_{j_3=1}^N e^{i \lambda_3\theta_{\ell_3}(t)}\right) e^{i\lambda_4\theta_m(t)}\right)\\
    \label{eq:HOKuramoto_fastRHS}
    &=\operatorname{Im}\left( r_{\lambda_1}~r_{\lambda_2}~r_{\lambda_3}~ e^{i\lambda_4\theta_m(t)} \right).
\end{align}
As in the section about the Kuramoto model \ref{sec:Kuramoto}, $r_m$ is the $m$-th order parameter
\begin{align*}
    r_m = \frac{1}{N} \sum_{\ell=1}^N e^{i\theta_\ell}.
\end{align*}
As can easily be seen precomputing the order parameters $r_m$ for $m=\lambda_1,\lambda_2,\lambda_3$ requires a $\mathcal O(N)$ function evaluations. Subsequently evaluating the right-hand side according to the formula \eqref{eq:HOKuramoto_fastRHS} takes another $\mathcal O(N)$ operations. Thus, in summary the complexity of this algorithm is linear in $N$, which is a significant reduction from the naive algorithm that scales with $N^4$. This shows the power of our approach: as long as one can exploit dense coupling structure, even higher-order (or polyadic, or hypergraph) interactions can be reduced from a high polynomial computational complexity in $N$ to linear complexity.

Even though there exist formulas that give a fast evaluation of the right-hand side only by using real numbers, deriving these formulas requires addition theorems on $\sin(\alpha_1+\alpha_2+\alpha_3+\alpha_4)$ and consequently these formulas tend to be very long, which is why we recommend the complex formula \eqref{eq:HOKuramoto_fastRHS}.

\subsection{Bornholdt-Rohlf Model}\label{sec:Bornholdt-Rohlf}

The Bornholdt-Rohlf model \cite{BornholdtRohlf} on a static all-to-all network is a discrete dynamical system with the iteration rule
\begin{align*}
    f_m(t) &= \sum_{\ell=1}^N v_\ell(t) + \mu v_m(t) + \sigma r_m, \quad r_m\sim \mathcal N(0,1),\\
    v_m(t+1) &= \operatorname{sgn}[f_m(t)],
\end{align*}
where $\sigma\ge 0$ is a parameter for the noise and $\mathcal N(0,1)$ denotes the standard normal distribution. Here, one needs to observe that the decisive sum in the definition of $f_m$ is independent of $m$. Therefore, this sum can be precomputed and then reused for each calculation of $f_m(t)$. In this way one can construct an algorithm whose complexity is linear in $N$. This example aims to illustrate in a simple setting that the neither the continuous-time assumption, nor the assumption about a particular ordinary differential equation structure matter. What does matter for being able to the CIA approach is that the computational bottleneck in a naive approach arises due to summing at each node over all its inputs.

\section{Extensions of the Community Structure Exploitation step (E1)}\label{sec:ExtensionE1}

The Community Structure Exploitation step (E1) as we have described it in Section \ref{sec:CIA} assumes that the nodes in each community are all-to-all coupled. However, there are cases for which the right-hand side of 
\begin{align*}
    x_m'(t) = \frac{1}{N} \sum_{\ell = 1}^N a_{m\ell}~ g(x_\ell(t), x_m(t))
\end{align*}
can be efficiently evaluated in step (E1) even though the graph represented by the adjacency matrix $A = (a_{m\ell})_{m,\ell}$ does not represent an all-to-all coupling or a vey dense coupling. Examples include rank one matrices and nearest-neighbor networks. In the following two subsections we briefly explain how an efficient evaluation on these networks is possible. For the sake of simplicity we assume that $g$ is of the form $g(x,y) = h(x-y)$ and $h$ consists of only one complex harmonic, i.e. $h(x) = e^{i x}$. The general case can then be obtained by approximating $h$ with more Fourier modes treating each harmonic individually.

\subsection{Rank One Coupling}

Here, we consider the case that the adjacency matrix is given by an outer product $a_{m\ell} = \alpha_m~\beta_\ell$ for two vectors $\alpha,\beta\in \R^N$, i.e. we focus on the system 
\begin{align}\label{eq:rank1system}
    x_m'(t) = \frac 1N \sum_{\ell = 1}^N \alpha_m\beta_\ell~ e^{i(x_\ell-x_m)}.
\end{align}
This is the only case in which we do not require $a_{m\ell}\in \{0,1\}$. Since the CIA evaluation step (E1) as described in Section \ref{sec:CIA} is a special case of a rank one coupling, i.e. when $\alpha_m=\beta_\ell=1$ for all $m,\ell$, this presents a generalization of (E1).
Based on \eqref{eq:rank1system}, we immediately see that it makes sense to precompute
\begin{align}\label{eq:rank1precompute}
    r = \frac{1}{N} \sum_{\ell = 1}^N \beta_\ell~e^{ix_\ell}.
\end{align}
Then, \eqref{eq:rank1system} simplifies to
\begin{align}\label{eq:rank1fast}
    x_m'(t) = \alpha_m~r~e^{-ix_m}.
\end{align}
Note that both \eqref{eq:rank1precompute} and \eqref{eq:rank1fast} can be evaluated with a complexity of $\mathcal O(N).$

\subsection{Nearest-Neighbor Coupling}
For a given $N\in \N$ and $k\in \N$ one can define a $k$-nearest-neighbor graph on $N$ nodes in terms of the adjacency matrix $A = (a_{m\ell})$ by setting
\begin{align*}
    a_{m\ell} = \begin{cases}1 &\quad \text{if } c_N(m,\ell)\le k\\ 0 &\quad \text{else}\end{cases},
\end{align*}
where $c_N$ is a circular distance given by
\begin{align*}
    c_N(m,\ell) = \min(\abs{m-\ell},~ N-\abs{m-\ell}).
\end{align*}
In other words, one can imagine all $N$ nodes uniformly distributed on a unit circle and then couple each node to its $k$ nearest neighbors in both directions. That means, we consider the system
\begin{align}
    \nonumber
    x_m'(t) &= \frac 1N \sum_{\ell = 1}^N a_{m\ell} e^{i(x_\ell-x_m)}\\
    \label{eq:nearestneighbor}
    &= \frac{1}{N}\sum_{\ell = m-k}^{m+k} e^{i(x_\ell-x_m)},
\end{align}
where the particle index in \eqref{eq:nearestneighbor} has to be understood modulo $N$. Since for $k\ge N/2$ this is an all-to-all coupling, the nearest-neighbor coupling is again a generalization of (E1) as described in Section \ref{sec:CIA}. However, the procedure to efficiently evaluate \eqref{eq:nearestneighbor} will be different and especially not based on precomputations but on an iterative method. In particular, given the representation \eqref{eq:nearestneighbor} one obtains
\begin{align*}
    x_m'(t) e^{ix_m} = \frac{1}{N} \sum_{\ell = m-k}^{m+k} e^{ix_\ell} =: F_m(x).
\end{align*}
With this notation
\begin{align}\label{eq:nearestneighboriteration}
    F_{m+1}(x) = F_m(x) -\frac 1N ( e^{ix_{m-k}} - e^{ix_{m+k+1}}),
\end{align}
which gives an iterative procedure to compute $F$. The steps to compute $x_m'(t)$ efficiently are therefore given by
\begin{enumerate}
    \item Compute $F_1(x)$ by the definition
    \begin{align*}
        F_1(x) = \frac 1N \sum_{\ell = 1-k}^{1+k} e^{ix_\ell}.
    \end{align*}
    \item Use formula \eqref{eq:nearestneighboriteration} to iteratively compute $F_2(x),F_3(x), \dots, F_N(x)$.
    \item Finally compute $x_m'(t) = F_m(x) e^{-ix_m}$.
\end{enumerate}
Note that the first step has a complexity of $\mathcal O(N)$. Each iteration in step 2 consists of only a finite number of operations $\mathcal O(1)$. Since $N-1$ iterations are necessary, the second step is of total complexity $\mathcal O(N)$. Finally the third step is obviously of complexity $\mathcal O(N)$ as well. Therefore, \eqref{eq:nearestneighbor} can be computed for all $m=1,\dots,N$ requiring a complexity of only $\mathcal O(N)$. Importantly, this is independent of $k$ and still holds when for example $k=rN$, where $r\in (0,1/2)$ is a factor that describes the coupling range. In this case the adjacency matrix $A$ has approximately $2rN^2$ non-zero entries, yet a efficient computation in $\mathcal O(N)$ is possible.

\end{document}